\documentclass{m2an}
\usepackage{amssymb}
\usepackage{float}
\usepackage{amsmath}
\usepackage{bbm}
\usepackage{multirow}	
\usepackage{enumerate}
\usepackage{graphicx}
\usepackage{graphics}
\usepackage{color}
\usepackage{subfigure}

\usepackage{amsthm}

\def\qin{{\quad\hbox{in}\quad}}
\def\qon{{\quad\hbox{on}\quad}}

\newcommand{\pin}{\mathbf{\Pi}}
\newcommand{\qy}{\quad\text{and}\quad}
\newcommand{\pr}{\partial}

\newcommand{\Scalk}{\Scal^K}

\newcommand{\tr}{\text{\rm tr}}

\newcommand{\kt}{{\widetilde{K}}}

\newcommand{\Om}{\Omega}
\newcommand{\G}{\Gamma}

\newcommand{\ds}{\displaystyle}
\newcommand{\alphan}{{\boldsymbol{\alpha}}}

\newcommand{\bkappa}{{\boldsymbol\kappa}}

\newcommand{\bsi}{{\boldsymbol\sigma}}

\newcommand{\btau}{{\boldsymbol\tau}}

\newcommand{\brho}{{\boldsymbol\rho}}
\newcommand{\bzeta}{{\boldsymbol\zeta}}

\newcommand{\bx}{{\boldsymbol{x}}}
\newcommand{\bcero}{{\boldsymbol{0}}}

\newcommand{\bt}{{\mathbf{t}}}

\newcommand\bH{\mathbf{H}}

\newcommand\rmC{\mathrm{C}}
\newcommand\rmH{\mathrm{H}}

\newcommand{\Prm}{\mathrm{P}}
\newcommand{\R}{\mathrm{R}}
\newcommand{\rmL}{\mathrm{L}}

\newcommand{\bbB}{\mathbb{B}}

\newcommand{\Acal}{{\mathcal{A}}}
\newcommand{\Bcal}{{\mathcal{B}}}
\newcommand{\Ecal}{{\mathcal{E}}}
\newcommand{\Ical}{{\mathcal{I}}}

\newcommand{\Gcal}{{\mathcal{G}}}

\newcommand{\Ocal}{{\mathcal{O}}}

\newcommand{\Scal}{\mathcal{S}}

\newcommand{\Tcal}{\mathcal{T}}

\newcommand{\bdiv}{\mathbf{div}}

\newcommand{\bcurl}{\mathbf{curl}}

\newcommand{\brot}{\mathbf{rot}}
\renewcommand{\div}{\mathrm{div}}
\newcommand{\rot}{\mathrm{rot}}

\newtheorem{assump}[thrm]{Assumption}

\begin{document}

\title{A posteriori error estimates for  mixed virtual element methods}\

\author{Andrea Cangiani}\address{Department of Mathematics, University of Nottingham, University Park, Nottingham NG7 2RD, UK,\\ email: {\tt Andrea.Cangiani@nottingham.ac.uk}}
\author{Mauricio Munar}\address{CI$^2$MA and Departamento de Ingenier\'\i a Matem\' atica,
                  Universidad de Concepci\' on, Casilla 160-C, Concepci\' on, Chile,\\
               email: {\tt  mmunar@ci2ma.udec.cl}}

\date{}

\begin{abstract}
We present an a posteriori error analysis for the mixed virtual element method (mixed VEM) applied to second order elliptic equations in divergence form with mixed boundary conditions.  The resulting error estimator is of residual-type. It only depends on quantities directly available from the VEM solution  and applies on very general polygonal meshes.
The proof of the upper bound relies on a global inf-sup condition, a suitable Helmholtz decomposition, and the local approximation properties of a Cl\'ement-type interpolant. In turn, standard inverse inequalities and localization techniques based on bubble functions are the main tools yielding the lower bound. Via the inclusion  of a fully local postprocessing of the mixed VEM solution,  we also show that the estimator provides a reliable and efficient control on the broken $\rmH(\div)$-norm error between the exact and the postprocessed flux.  
Numerical examples  confirm the theoretical properties of our estimator, and show that it can be effectively used to drive an adaptive mesh refinement algorithm.  
\end{abstract}


\maketitle

\section{Introduction}

The Virtual Element Method (VEM) was originally introduced in \cite{bbcmmr-2013-MMMA} for the solution of elliptic problems, followed by the mixed VEM proposed in~\cite{bfm-2014-ESAIM}. 
Subsequently, new mixed VEMs have been analysed for the  solution of the Stokes, Navier-Stokes, and Brinkman problem~\cite{pbmv-2014-SIAM,cvm-2016-SIAM,cg-2017-IMA,blv-2017-ESAIM,blv-2018-SIAM,gms-2018-M3AS,caceres_brinkman,v-2018-M3AS}.

One of the defining characteristics of the VEM is that it allows for the use of very general polygonal and polyhedral meshes. As such, it naturally lends itself as a flexible solution step within automatically adaptive algorithms.
Indeed, mesh refinement and coarsening strategies can  be  implemented  very  easily and efficiently as, for instance,  hanging nodes are simply treated as new nodes, with no detrimental affect on the quality of the approximation.
Moreover, it has been shown that the VEM in primal form allows for extremely aggressive mesh adaptation producing strongly solution-adapted polygonal meshes~\cite{cgto-2017-NM}. 

In this respect, the design and analysis of adaptive mesh refinement strategies based on robust a posteriori error indicators for the VEM approach and, in particular, for the mixed-VEM is an attractive proposition.

Several error estimators have been proposed in the context of VEM for primal forms (see,  e.g.,\cite{bm-2015-ESAIM,cgto-2017-NM,bb-2017-MMMA,mrr-2017-CMA,mr-2018-IMA,cbg-2019-CMAM}).  Firstly, the authors of \cite{bm-2015-ESAIM} proposed a posteriori error bounds for the $C^1$-conforming VEM for the two-dimensional Poisson problem. Next, a posteriori error bounds for the $C^0$-conforming VEM for the discretization of second-order linear elliptic reaction-convection-diffusion problems with nonconstant coefficients in two and three dimension were proposed in \cite{cgto-2017-NM}, whereas a residual-based a posteriori error estimator for the VEM discretization of the Poisson problem with discontinuous diffusivity coefficient was introduced and analysed in \cite{bb-2017-MMMA}. Moreover, in \cite{mrr-2017-CMA} and \cite{mr-2018-IMA}, the authors developed a  posteriori error analysis of a VEM approach for the Steklov eigenvalue problem and the spectral analysis  for the elasticity equations, respectively. Finally, in \cite{cbg-2019-CMAM} a general recovery-based a posteriori error estimation framework for the VEM of arbitrary order on general polygonal/polyhedral meshes has been developed. A posteriori analises of other techniques of mixed-type on general meshes have been presented in~\cite{DIPIETRO2016} for the Mixed High-Order method, in~\cite{MFD-apost} for the Mimetic Finite Difference method, and in ~\cite{vohralik-2008} for lowest-order locally conservative methods. However, to the best of our knowledge, \emph{no} a posteriori error analysis for mixed VEM is available from the literature.

The  aim  of  this  paper  is  to  introduce the basic tools to develop the a posteriori error analysis for the mixed VEM. To this end, we consider a second order elliptic equation in divergence form with mixed boundary condition, discretised using  the basic  mixed VEM of~\cite{bbmr-2016-ESAIM}. 
As usual in the VEM approach, we introduce fully computable approximations for the virtual approximation of the  flux variable and establish its corresponding a priori error estimates. In particular, in order to improve on the sub-optimal order provided by the computable component of the flux variable in the broken $\rmH(\div)$-norm,  observed numerically in \cite{gms-2018-c}, we follow \cite{caceres_brinkman} and construct by postprocessing  second computable approximation of the flux variable, which has an optimal rate of convergence in the  aforementioned norm.

The a posteriori error analysis is based on a global inf-sup condition coming from the well-posedness of the continuos problem. Upper bounds are shown for the scalar variable in the $\rmL^2$-norm, the VEM flux variable in the $\rmH(\div)$-norm, its projection in the $\rmL^2$-norm, and postprocessing in the broken $\rmH(\div)$-norm. The proof   uses  properties  of the interpolation operator associated to the virtual subspace of the flux variable and Cl\'ement-type interpolation operators, together with a suitable Helmholtz decomposition.  Moreover,  some  inverse  inequalities and localization techniques based on bubble functions will serve to show a lower bound for the error. In this way, we are able to establish the equivalence up to virtual inconsistency terms between the  error and the error estimator for the postprocessing of the  virtual element approximation, measured in the broken $\rmH(\div)$-norm.

\subsection{Outline}

The remainder of the paper has been structured as follows. In what is left of this section,
we introduce some standard notations and the required functional spaces. In Section~\ref{sec:modelprob} we introduce the model problem and presents the associate variational formulation.  In Section~\ref{sec:method}, we present the  mixed virtual element scheme. The a posteriori error analysis is laid down in details
in Section~\ref{sec:aposteriorianalysis}. In Section ~\ref{sec:numerics}, we propose an adaptive algorithm and test its effectiveness with some numerical
examples. Finally, in Section~\ref{sec:conclusions} we give some concluding remarks.

\subsection{Preliminaries}

Let us assume that $\Om\subset\R^2$ be a bounded domain with polygonal boundary $\G$. We denote by $\nu$ the outward unit normal vector to the boundary $\G$. Moreover, we assume that $\G$ admits a disjoint partition $\G=\overline{\G}_D\cup\overline{\G}_N$, where $\G_D$ and $\G_N$ are open subsets of $\G$, with $|\G_D|,|\G_N|\neq 0$.

For $s\geq 0$, the symbol $|\cdot|_{s,\Omega}$
stands for the norm of the Hilbertian Sobolev spaces $\rmH^s(\Om)$,  with the convention $\rmH^0(\Omega):=\rmL^2(\Om)$. We also define the
Hilbert space
\begin{equation*}
\rmH(\div;\Om):=\Big\{\btau\in[\rmL^2(\Om)]^2:\quad\div\,\btau\in\rmL^2(\Om)\Big\}\,,
\end{equation*}
whose norm is given by $\|\btau\|_{\div;\Om}^2:=\|\btau\|_{0,\Om}^2+\|\div\,\btau\|_{0,\Om}^2$. Hereafter, we use the following notation for any vector field $\btau=(\tau_i)_{i=1,2}$ and any scalar field $v$:
\begin{equation*}
\div\,\btau:= \partial_1\tau_1 + \partial_2\tau_2\,\quad\rot\,\btau :=\partial_1\tau_2 - \partial_2\tau_1\qy\brot\,v:=\left(\partial_2 v,-\partial_1 v\right)^\bt\,.
\end{equation*}
Additionally, we need to introduce the following spaces
\begin{equation}\label{space-hgn}
H:=\Big\{\btau\in\rmH(\div;\Om):\quad\btau
\cdot\nu = 0 \qon\G_N\Big\}\qy Q:=\rmL^2(\Om)\,,
\end{equation}
endorsed with the norms
\begin{equation*}
\|\tau\|_H:=\|\tau\|_{0,\Om} + \|\div\,\tau\|_{0,\Om}\qy\|v\|_Q:=\|v\|_{0,\Om}\,.
\end{equation*}
Furthermore, we make use of the product space $H\times Q$ with the norm
\begin{equation*}
\|(\tau,v)\|_{H\times Q}:=\|\tau\|_H + \|v\|_Q\,.
\end{equation*}

In addition, we will denote with $c$ and $C$, with or without subscripts, tildes, or hats, a generic constant
independent of the mesh parameter $h$, which may take different values in different occurrences.

\section{The model problem}\label{sec:modelprob}

We consider the problem
\begin{equation}\label{problem}
-\div(\bkappa\nabla u)=f\qin\Om,\quad u=g\qon\G_D\qy (\bkappa\nabla u)\cdot\nu=0 \qon\G_N,
\end{equation}
where $f\in \rmL^2(\Om),\,g\in\rmH^{1/2}(\G_D)$ and $\bkappa\in[\rmL^{\infty}(\Om)]^{2\times 2}$ is an uniformly positive definite tensor, which is assumed to be known. In particular, we denote by $\kappa^{*}$
the positive constant satisfying
\begin{equation}
\bkappa^{-1}\boldsymbol{\zeta}\cdot\boldsymbol{\zeta}\geq \kappa^{*}|\boldsymbol{\zeta}|^2\,,\quad\,\forall\,\boldsymbol{\zeta}\in[\rmL^2(\Om)]^2\,.
\end{equation}

By introducing the flux variable $\bsi:=\bkappa\nabla u$ in $\Om$ as  additional unknown, a mixed variational formulation of \eqref{problem} becomes:

Find $(\bsi,u)\in H\times Q$ such that
\begin{equation}\label{scheme_c}
\begin{array}{rrcl}
a(\bsi,\btau)+b(\btau,u)&=&\langle\btau\cdot\nu,g\rangle_{\G_D}&\quad\forall\;\btau\in H,\\[2 ex]   
b(\bsi,v)&=&-\ds\int_\Om fv&\quad\forall\;v\in  Q\,,
\end{array}
\end{equation}
where $\langle\cdot,\cdot\rangle$ stands for the duality pairing between $\rmH^{-1/2}(\G_D)\to\rmH^{1/2}(\G_D)$. In turn, $a:H\times H\to\R$ and $b:H\times Q\to\R$
 are the bounded bilinear forms defined by
\begin{equation}\label{forms}
a(\bsi,\btau):=\ds\int_\Om\bkappa^{-1}\bsi\cdot\btau,\qy
b(\btau,u):=\ds\int_\Om u\,\div\,\btau.
\end{equation}

Under the assumptions on $\bkappa,f$ and $g$, the  existence and uniqueness
of the weak solution of \eqref{scheme_c} is consequence of the Bab\v{u}ska-Brezzi theory.

\section{The virtual element method}\label{sec:method}

Let $\{\Tcal_h\}_{h>0}$ be a  family of decompositions of $\Om$ into open non-overlapping  polygonal elements. Then, for each $K\in\Tcal_h$ we denote its diameter by $h_K$, and also, as usual, $h:=\max\Big\{h_K : K\in\Tcal_h\Big\}$. In what follows we make the following mesh regularity assumptions which are standard in this context (cf.\cite{bbcmmr-2013-MMMA,bfm-2014-ESAIM}).

\begin{assump}\label{ass:mesh}
The family of decompositions $\lbrace\Tcal_h\rbrace_{h>0}$ satisfies:
	\begin{enumerate}
	\item[a)] the ratio between the shortest edge and the diameter $h_K$ of $K$ is bigger
	than $C_\Tcal$, and
	\item[b)] $K$ is star-shaped with respect to a ball $B$ of radius $C_\Tcal h_K$ and
	center $\bx_B\in K$.
	\end{enumerate}
\end{assump}

\begin{rmrk}\label{remark:meshassumptions}
The above assumptions imply that each $K\in\Tcal_h$ is simply connected and that there exists an integer $N_\Tcal$ (depending only on $C_\Tcal$), such that the numbers of edges of each $K\in\Tcal_h$ is bounded above by $N_\Tcal$.

Moreover, as each element $K$ is star-shaped, it admits a sub-triangulation $\Tcal_h^K$  obtained by joining each vertex of $K$ with a point with respect to which $K$ is starred. And the uniform bound on the diamater of the mesh edges ensures that the resulting global triangulation $\widehat{\Tcal_h}:=\ds\bigcup_{K\in\Tcal_h}\Tcal_h^K$ is shape-regular.

We finally note that the  above assumptions allow for very general possibly non-convex polygonal elements. In particular, they permit the natural incorporation of so-called hanging nodes, thus completely avoiding the need of removing hanging nodes typical of standard mesh adaptation algorithms.
\end{rmrk}

Now, given an integer $\ell\geq 0$ and  $\Ocal\subseteq\R^d$, $d=1,2$, we denote by $\Prm_\ell(\Ocal)$ the space of polynomials on $\Ocal$ of degree up to $\ell$. Then,  given  an edge $e\in\pr K$ with barycentric $x_e$ and diameter $h_e$, we denote
the following set of $(\ell+1)$ normalized monomials on $e$
\begin{equation*}
 \Bcal_\ell(e)\ :=\ \left\{\left(\frac{x-x_e}{h_e}\right)^j\right\}_{0\leq j\leq \ell}\,,
 \end{equation*}
which certainly constitutes a basis on $\Prm_{\ell}(e)$. Similarly, on $K\in\Tcal_h$ with barycenter $\bx_K$, we define the following
set of $\frac{1}{2}(\ell+1)(\ell+2)$ normalized monomials
\begin{equation*}
\Bcal_\ell(K)\ :=\ \left\{\left(\frac{\bx-\bx_K}{h_K}\right)^{\alphan}\right\}_{0\leq|\alphan|\leq \ell}\,,
\end{equation*}
which is a basis of $\Prm_{\ell}(K)$. Notice that in the definition of $\Bcal_\ell(K)$ above,
we made use of the multi-index notation, that is, given $\bx:=(x_1,x_2)^{\bt}\in \R^2$
and $\alphan:=(\alpha_1,\alpha_2)^{\bt}$, with non-negative integers $\alpha_1,\alpha_2$, we
set $\bx^{\alphan}:=x_1^{\alpha_1}x_2^{\alpha_2}$ and $|\alphan|:=\alpha_1+\alpha_2$.

We further let $\Gcal_{\ell}(K)$ be a basis of
$\big(\nabla\Prm_{\ell+1}(K)\big)\cap[\Prm_{\ell}(K)]^{2}$, whereas with $\Gcal_{\ell}^\perp(K)$ we denote a basis  of the $[\rmL^{2}(K)]^{2}$-orthogonal
of $\Gcal_{\ell}(K)$ in $[\Prm_{\ell}(K)]^{2}$.

Throughout the paper, we denote by
 $\Pi_k^0:\rmL^2(K)\to\Prm_k(K)$ the $\rmL^2(K)$-orthogonal projection onto the space $\Prm_k(K)$, for any $K\in\Tcal_h$ and $k\geq 0$. In addition, we will make use of a vectorial version of the aforementioned projector, which is denoted by $\pin_k^0$. The following approximation properties of these projectors are well-known:
 \begin{equation}\label{approx-pk}
\Vert v-\Pi_k^0(v)\Vert_{0,K}\ \leq\ Ch_K^m|v|_{m,K}\qy \Vert \btau-\pin_k^0(\btau)\Vert_{0,K}\ \leq\ Ch_K^m|\btau|_{m,K}
\end{equation}
for all $K\in\Tcal_h$, and for all $v\in\rmH^m(K)$, $\btau\in[\rmH^{m}(K)]^{2}$, with $m\in\lbrace 0,1,\dots,k+1\rbrace$.

\section{Virtual subspaces and its approximation properties}\label{proper}

For any integer $k\geq 0$, we introduce the finite dimensional subspaces of $H$ and $Q$, respectively, given by
\begin{equation}\label{space-hh}
H_h\ :=\ \Big\{\btau\in H:\quad\btau\big\vert_K\in H_h^K\quad\forall\;K\in\Tcal_h\Big\},
\end{equation}
and
\begin{equation}\label{space-qh}
Q_h\ :=\ \Big\{ v\in Q:\quad v\big\vert_K\in Q_h^K\quad\forall\;K\in\Tcal_h\Big\},
\end{equation}
where  $Q_h^K:=\Prm_{k}(K)$, and $H_h^K$ is the virtual element space introduced in \cite[Section 3.1]{bbmr-2016-ESAIM}. This is defined by
\begin{equation}\label{space-hk}
\begin{array}{rcl}
H_h^{K} & := & \Big\{\btau\in\rmH(\div;K)\cap\rmH(\rot;K)\,:\quad \btau\cdot\nu|_e\in\Prm_k(e)\quad\forall\ \text{edge }e\in\pr K\,,\quad\\[1ex]
& & \phantom{\btau\in\rmH(\div;K)\,:} \quad\div\,\btau\in\Prm_{k}(K)\qy\rot\,\btau\in\Prm_{k-1}(K)\Big\}.
\end{array}
\end{equation}
and is characterised by the following degrees of freedom (cf.~\cite{bbmr-2016-ESAIM,bbmr-NM-2016}):
\begin{equation}\label{dof-hh}
\begin{array}{lll}
\displaystyle\int_e q\,(\btau\cdot\nu)&\;\; \forall\ q\in\Bcal_k(e)\,, & \forall\ \mbox{edge}\ e\qin\Tcal_h,\\[3ex]
  \displaystyle\int_K\btau\cdot\nabla q & \forall\; q \in\Bcal_{k}
  (K)\setminus\lbrace 1\rbrace,&\forall\;K\in\Tcal_h,\\[3ex]
 \displaystyle\int_K\btau\cdot \eta& \forall\ \eta\in\Gcal^{\perp}_{k}
  (K)\,, &\forall\;K\in\Tcal_h.
\end{array}
\end{equation}

As was remarked in \cite[Section 3.2]{bbmr-2016-ESAIM} (see also \cite[Section 3.5]{bbmr-NM-2016}), the degrees of freedom \eqref{dof-hh} allow the explicit computation of the projection $\pin_k^0(\btau)$ using only the degrees of freedom  of $\btau$.
Moreover, collected together, the local degrees of freedom~\eqref{dof-hh} provide a set of degrees of freedom for the global virtual element space $H_h$.
  
For each $\btau\in H$ such that $\btau\big\vert_K\in[\rmH^1(K)]^{2}$ for all $K\in\Tcal_h$, we may denote by  $\btau_I\in H_h$ the Lagrange interpolant of $\btau$ with respect to the  degrees of freedom~\eqref{dof-hh}.
For each $q\in\Bcal_{k}(K)$ we find that
\begin{equation*}
\ds\int_K q\,\div(\btau-\btau_I)\ =\ -\ds\int_K(\btau-\btau_I)\cdot\nabla q+\ds\int_{\pr K}q(\btau-\btau_I)\cdot\nu\ =\ 0,
\end{equation*}
which, thanks to the fact that $\div\,\btau_I\in\Prm_{k}(K)$, implies the commutative property
\begin{equation}\label{conmutativo}
\div\,\btau_I\ =\ \Pi_k^0(\div\,\btau)\qquad\forall\ \btau\in[\rmH^1(K)]^{2}.
\end{equation}
Hence we have the following approximation error estimates~\cite{bbcmmr-2013-MMMA,bbmr-2016-ESAIM}.
\begin{lmm}\label{approx-interpol-Hh}
Let $r$ be an integer such that 
$1\leq r\leq k+1$. Then, there exists a constant $C>0$,
independent of $K$, such that for each $\btau\in[\rmH^r(K)]^{2}$ such that
$\div\,\btau\in\rmH^r(K)$ there holds
\begin{equation}\label{approx-hdiv}
\|\btau -\btau_I\|_{\div;K}\ \leq\ C\,h_K^r\,\Big\{|\btau|_{r,K}\, +\, |\div\,\btau|_{r,K}\Big\}\qquad\forall\ K\in\Tcal_h\,.
\end{equation}
\end{lmm}
\begin{proof}
The bound on the divergence term follows from~\eqref{conmutativo}  and~\eqref{approx-pk}. The result then follows from classical arguments~\cite{brenner}. 
\end{proof}

\subsection{Discrete formulation}\label{form-disc}

We now aim to define a virtual scheme for our problem (\ref{scheme_c}) based on the discrete spaces \eqref{space-hh} and \eqref{space-qh}. To this end, we first notice that the bilinear form $b$ (cf.\eqref{forms}) is explicitly computable for all $(\btau,v)\in H_h\times Q_h$, just by accessing the degrees of freedom \eqref{dof-hh}. On the contrary, for each $K\in\Tcal_h$, the local version $a^K:H_h^K\times H_h^K \to\R$ of the bilinear form $a$,  which, is defined for all $\bzeta,\btau\in H_h^K\times H_h^K$ by
\begin{equation}\label{forma-ak}
a^K(\bzeta,\btau):=\ds\int_K\bkappa^{-1}\bzeta\cdot\btau,
\end{equation}
is not explicitly computable for $\bzeta,\btau\in H_h^K$ since in general $\bzeta$ and $\btau$ are not known explicitly on the whole of $K$. In order to deal with this difficulty, we follow \cite[Section 3.3]{bbmr-2016-ESAIM} and  introduce a local bilinear form 
 $a_h^K:H_h^K\times H_h^K \to\R$ defined by
\begin{equation}\label{forma-ahk}
a_h^K(\bzeta,\btau):=a^K(\pin_k^0(\bzeta),\pin_k^0(\btau))+\Scalk(\bzeta-\pin_k^0(\bzeta),\btau-\pin_k^0(\btau))\,,
\end{equation}
where 
$\Scal^K:H_h^K\times H_h^K\to\R$ is any symmetric and positive  definite bilinear form such that
\begin{equation}\label{cota-scalk}
\widehat{c}_0a^K(\bzeta,\bzeta) \leq \Scal^K(\bzeta,\bzeta)\ \leq\ \widehat{c}_1a^K(\bzeta,\bzeta)\qquad\forall\;\bzeta\in H_h^K,\quad\mbox{with}\quad\pin_k^0(\bzeta)=0\,,
\end{equation}
with constants $\widehat{c}_0,\widehat{c}_1>0$ which depend only on the shape regularity constant $C_\Tcal$ and on $\bkappa$.  In particular, to define $\Scal^K$ we can consider the bilinear form associated to the identity matrix in $\R^{n_k^K}$ with respect to the local basis determined by the degrees of freedom \eqref{dof-hh}, and where $n_k^K=\dim H_h^K$. (cf.\cite{bbcmmr-2013-MMMA,bfm-2014-ESAIM})

The following two lemmas establish the properties of the bilinear form $a_h^K$ and the consistency error between $a_h^K$ and $a^K$, respectively.
\begin{lmm}\label{prope-ah}
For all $K\in\Tcal_h$, there holds
\begin{equation*}
{\bf{(Consistency)}}\qquad a_h^K(p,\bzeta)=\ds\int_K\bkappa^{-1}p\cdot\pin_k^0(\bzeta)\qquad\forall\ p\in[\Prm_k(K)]^2\qy\forall\ \bzeta\in H_k^K\,,
\end{equation*}
and further, there exist constants $\alpha_{*},\alpha^{*}>0$, such that
\begin{equation*}
{\bf{(Stability)}}\qquad\alpha_{*} a^K(\bzeta,\bzeta)\leq  a_h^K(\bzeta,\bzeta)\leq \alpha^{*} a^K(\bzeta,\bzeta)\qquad\forall\ \bzeta\in H_h^K,\forall\ K\in\Tcal_h\,.
\end{equation*}
\end{lmm}
\begin{proof}
We refer to \cite{bbmr-2016-ESAIM} and \cite{bm-2014-IMA} for the details.
\end{proof}
\begin{lmm}\label{inconsistency}
There exists a constant $C>0$, depending only on $\bkappa\,,\widehat{c_1}$ and $\alpha^{*}$, such that
\begin{equation*}
(a_h^K-a^K)(\bzeta,\btau)\leq C\Big\{\Vert\bzeta-\pin_k^0(\bzeta)\Vert_{0,K}+\|\bkappa^{-1}\pin_k^0(\bzeta)-\pin_k^0(\bkappa^{-1}\pin_k^0(\bzeta))\|_{0,K}\Big\}\Vert\btau\Vert_{0,K}
\end{equation*}
for all $\bzeta,\btau\in H_h^K$ and for all $K\in\Tcal_h$.
\end{lmm}
\begin{proof}
We have that
\begin{equation*}
\begin{array}{lll}
(a_h^K-a^K)(\bzeta,\btau)&=&-\ds\int_K\Big\{\bkappa^{-1}\pin_k^0(\bzeta)-\pin_k^0(\bkappa^{-1}\pin_k^0(\bzeta))\Big\}\cdot(\btau-\pin_k^0(\btau))
-\ds\int_K\bkappa^{-1}(\bzeta-\pin_k^0(\bzeta))\cdot\btau\\[2 ex]
&&+\Scalk(\bzeta-\pin_k^0(\bzeta),\btau-\pin_k^0(\btau))\,.
\end{array}
\end{equation*}
The results now follows from Cauchy-Schwarz inequality and the properties of the bilinear $\Scalk$.
\end{proof}
According to the definition \eqref{forma-ahk} the global discrete bilinear form $a_h:H_h\times H_h\to\R$ can now be defined summing together the local contribution \eqref{forma-ahk}, that is
\begin{equation}\label{forma-ah}
a_h(\bzeta,\btau):=\sum_{K\in\Tcal_h}a_h^K(\bzeta,\btau)\qquad\forall\;\bzeta,\btau\in H_h.
\end{equation}
In this way, the virtual element method associated with the formulation (\ref{scheme_c}) reads: 

\noindent 
Find $(\bsi_h,u_h)\in H_h\times Q_h$ such that
\begin{equation}\label{scheme_d}
\begin{array}{rrcl}
a_h(\bsi_h,\btau_h)+b(\btau_h,u_h)&=&\langle\btau_h\cdot\nu,g\rangle_{\G_D}&\quad\forall\;\btau_h\in H_h, \\[2 ex]
              b(\bsi_h,v_h)&=&-\ds\int_\Om fv_h&\quad\forall\;v_h\in Q_h.
\end{array}
\end{equation}
\noindent 
The well-posedness of \eqref{scheme_d} follows from Lemma  \ref{prope-ah} and of the well-posedness of \eqref{scheme_c}. In addition, we have the following result about the a priori error estimates for the schemes \eqref{scheme_c} and \eqref{scheme_d}.

\begin{thrm}\label{apriori-1}
Let $(\bsi,u)\in H\times Q$ and $(\bsi_h,u_h)\in H_h\times Q_h$ be the unique solutions of the continuous  and discrete schemes \eqref{scheme_c} and \eqref{scheme_d}, respectively. In addition,  assume that for some $s\in [1, k+1]$ there hold $\bsi\big\vert_K\in\bH^s(K)$ and $\div\,\bsi\big\vert_K\,,u\big\vert_K\in\rmH^s(K)$ for  each $K\in\Tcal_h$. Then, there exist a positive constant $C>0$, independent of $h$, such that 
\begin{equation}
\Vert(\bsi,u)-(\bsi_h,u_h)\Vert_{H\times Q}\leq Ch^s\left\lbrace\sum_{K\in\Tcal_h}|\bsi|_{s,K}^2\, +\, |\div\,\bsi|_{s,K}^2+|u|_{s,K}^2\right\rbrace^{1/2}\,.
\end{equation} 
\end{thrm}
\begin{proof} 
The result is consequence of \cite[Theorem 6.1]{bfm-2014-ESAIM} and of a straightforward application of the approximation properties provided by \eqref{approx-pk} and Lemma~\ref{approx-interpol-Hh}.
\end{proof}

\subsection{Computable approximations}\label{sec:computable}

A first fully computable approximation $\widehat{\bsi}_h\in Q$ of the VEM solution $\bsi_h\in H$ is given by
\begin{equation}\label{approx-flux}
\widehat{\bsi}_h:=\pin_k^0(\bsi_h).
\end{equation}
The corresponding a  priori error estimates for the error $\|\bsi-\widehat{\bsi}_h\|_Q$ immediately follows from the foregoing Theorem~\ref{apriori-1} and the triangle inequality. 
\begin{thrm}\label{apriori-2}
Let $(\bsi,u)\in H\times Q$ and $(\bsi_h,u_h)\in H_h\times Q_h$ be the unique solutions of the continuous  and discrete schemes \eqref{scheme_c} and \eqref{scheme_d}, respectively. In addition,  assume that for some $s\in [1, k+1]$ there hold $\bsi\big\vert_K\in\bH^s(K)$ and $u\big\vert_K\in\rmH^s(K)$ for  each $K\in\Tcal_h$. Then, there exists a positive constant $C>0$, independent of $h$, such that 
\begin{equation}
\Vert\bsi-\widehat{\bsi}_h\Vert_{Q}+\Vert u-u_h\Vert_{Q}\leq   Ch^s\left\lbrace\sum_{K\in\Tcal_h}|\bsi|_{s,K}^2\,+|u|_{s,K}^2\right\rbrace^{1/2}\,.
\end{equation}
\end{thrm}

Next, motivated by the non-satisfactory order
provided by $\widehat{\bsi}_h$ in the broken $\rmH(\div)$-norm (see \cite[Section 5] {gms-2018-c} for numerical evidences of this fact), we proceed as in \cite[Section 5.3] {caceres_brinkman} (see also~\cite{cgs-2018-SIAM}) and construct, by local postprocessing, a second approximation $\bsi_h^{\star}$ for the flux variable $\bsi$ which has an optimal rate of convergence in such norm. To this end, for each $K\in\Tcal_h$ we let $(\cdot,\cdot)_{\div;K}$ be the usual $\rmH(\div;K)$-inner product with induced norm $\|\cdot\|_{\div;K}$ and let $\bsi_h^{\star}\big\vert_K:=\bsi_{h,K}^{\star}\in[\Prm_{k+1}(K)]^2$ 
be the unique solution of the local problem
\begin{equation}\label{computo-smn-broken}
(\bsi_{h,K}^{\star},\btau_h)_{\div;K}=\ds\int_K\widehat{\bsi}_h\cdot\btau_h\ +\ \int_K\div\,\bsi_h\,\div\,\btau_h\qquad\forall\;\btau_h\in[\Prm_{k+1}(K)]^2.
\end{equation}
We stress that $\bsi_{h,K}^{\star}$ can be explicitly computed for each $K\in\Tcal_h$, independently.  Then, the rate of convergence for the broken $\rmH(\div;\Om)$-norm of $\bsi-\bsi_h^{\star}$ is established as follows.
\begin{thrm}\label{rate-broken}
Assume that the hypotheses of Theorem \Rref{apriori-1} are satisfied. 
Then, there exists a positive constant $C$, independent of $h$, such that
\begin{equation}\label{smn-broken}
\left\lbrace\sum_{K\in\Tcal_h}\|\bsi-\bsi_{h,K}^{\star}\|_{\div;K}^2\right\rbrace^{1/2}
\,\le\, Ch^s\left\lbrace\sum_{K\in\Tcal_h}|\bsi|_{s,K}^2\, +\, |\div\,\bsi|_{s,K}^2\right\rbrace^{1/2}\,.
\end{equation}
\end{thrm}
\begin{proof}
See \cite[Section 5.3,Theorem 5.5] {caceres_brinkman}.
\end{proof}

\section{A posteriori error analysis}\label{sec:aposteriorianalysis}

In this section we develop a residual-based a posteriori error analysis for the mixed virtual element scheme \eqref{scheme_d}.  The proof of the a posteriori upper bound on the error is based on a global inf-sup condition, (cf.~\cite{mgr-2016-MMNA}),  and a suitable Helmholtz decomposition; the lower bound is derived as usual via techniques based on bubble functions together with inverse inequalities.

\subsection{Preliminaries}\label{preli}

We let $\Ecal_h=\Ecal_h(\Om)\cup\Ecal_h(\G_D)\cup\Ecal_h(\G_N)$ be the set of all edges of $\Tcal_h$, where  $\Ecal_h(\Om):=\left\lbrace e\in\Ecal_h:\quad e\subseteq\Om\right\rbrace$, $\Ecal_h(\G_D):=\left\lbrace e\in\Ecal_h:\quad e\subseteq\G_D\right\rbrace$, and $\Ecal_h(\G_N):=\left\lbrace e\in\Ecal_h:\quad e\subseteq\G_N\right\rbrace$. And, for a given $K\in\Tcal_h$, we denote by $\Ecal(K)\subset\Ecal_h$ the set of edges of $K$.  
Given an edge $e\in\Ecal_h$, we let $h_e$ be its  length and we fix a unit normal vector $\nu_e:=(\nu_1,\nu_2)^\bt$ and let $s_e:=(-\nu_2,\nu_1)^\bt$ be the corresponding unit tangential vector along $e$. However, when no confusion arises, we simply write $\nu$ and $s$ instead of $\nu_e$ and $s_e$, respectively. Now, given $\bzeta\in[\rmL^2(\Om)]^{2}$, for each  $K\in\Tcal_h$ and $e\in\Ecal_h(\Om)\cap\Ecal(K)$ we denote  by $[\![\bzeta\cdot s]\!]$ the tangential jump of $\bzeta$ across $e$, that is $[\![\bzeta\cdot s]\!]:=(\bzeta\big\vert_K-\bzeta\big\vert_{K^\prime})\big\vert_e\cdot s$, where $K$ and $K^\prime$ are the elements of $\Tcal_h$ having $e$ as a common edge.

We first recall the conforming VEM spaces from \cite{bbcmmr-2013-MMMA}, which will be used as an auxiliary space in the a posteriori analysis below. Given $k\geq 0$, we consider the space defined by
\begin{equation*}
V_h:=\left\lbrace v\in \rmH^1(\Om):\quad v\big\vert_{\pr K}\in \bbB_{k}(\pr K)\mbox{ and }\Delta v\in \Prm_{k-1}(K)\quad\forall K\in\Tcal_h\right\rbrace,
\end{equation*}
where
\begin{equation*}
\bbB_{k}(\pr K)\ :=\ \big\{v\in \mathrm{C}(\pr K)\,:\quad v|_e\in \Prm_{k+1}(e)\quad\forall\ \text{edge }e\subseteq\pr K\big\}.
\end{equation*}
It has been shown in~\cite[Section 4, Proposition 4.2]{gonzalo} that there exists an interpolation operator $\Ical_h:\rmH^1(\Om)\to V_h$, such that there holds
\begin{equation}\label{ical-1}
\Vert v-\Ical_h(v)\Vert_{0,K}+h_K|v-\Ical_h(v)|_{1,K}\leq c_1h_K \Vert v\Vert_{1,K}\quad\forall\; v\in H^1(K).
\end{equation}
From this, using a scaled trace inequality, the Cauchy-Schwarz inequality, and Assumption~\ref{ass:mesh} it follows that 
\begin{equation}\label{ical-2}
\Vert v-\Ical_h(v)\Vert_{0,e}\leq c_2h_e^{1/2}\Vert v\Vert_{1,K}\quad\forall\;e\in\Ecal_h.
\end{equation}

We now let $\rmH_{\G_N}^1(\Om):=\left\lbrace v\in\rmH^1(\Om):\,
v=0\,\text{ on }\,\G_N\right\rbrace$ and consider the virtual element subspace given by 
\begin{equation}\label{tildeVh}
\widetilde V_{h}:=V_h\cap\rmH_{\G_N}^1(\Om).
\end{equation} 
Also, we introduce, analogously as before, the interpolation operator $\widetilde I_{h}:\rmH_{\G_N}^1(\Om)\to \widetilde V_{h}$ such that  $\widetilde I_{h}:=I_{h}\big\vert_{\rmH_{\G_N}^1(\Om)}$. In addition, the following lemma establishes an important relation between the virtual spaces $\widetilde V_{h}$ and $H_h$ (cf.\eqref{space-hh}).

\begin{lmm}\label{curl-hk}
For $k\geq 0$, given $v\in \widetilde V_{h}$ we have $\brot\,v\in H_h$.
\end{lmm}
\begin{proof}
Given $v\in \widetilde V_{h}$, it is easy to see that $\brot\,v\in H$. Moreover, given $K\in\Tcal_h$, we observe that $\rot(\brot\,v)=-\Delta v\in \Prm_{k-1}(K)$.
Furthermore, following \cite[Section 8, Theorem 3]{bbmr-NM-2016}, we have that $\brot\,v\cdot\nu\big\vert_e=\nabla v\cdot s\big\vert_e\in \Prm_{k}(e)$ for all edge $e\in\pr K$ . Hence, we conclude that $\brot\,v\big\vert_K\in H_h^K$ for all $K\in\Tcal_h$.
\end{proof}
We now recall from \cite[Section 3.3]{caceres_stokes} some preliminary notations and technical 
results. For each element $K\in\Tcal_h$ we first define $\kt \,:=\, T_K(K)$, where $T_K:\R^2\to \R^2$ is the 
bijective affine mapping defined by 
\[
T_K(\bx) \,:=\, \frac{\bx-\bx_B}{h_K}\quad\forall\;\bx\in \R^2\,.
\]
Then, as it was remarked in \cite[Section 3.3]{caceres_stokes}, it is easy to see that the diameter
$h_{\kt}$ of $\kt$ is $1$, the shortest edge of $\kt$ is bigger than $C_{\Tcal}$ (which follows from Assumption~\ref{ass:mesh}), and $\kt$ is star-shaped with  respect to a ball $\widetilde{B}$ of radius $C_{\Tcal}$ and centered at the origin. 
Then, by connecting each vertex of $\kt$ to the center of $\widetilde{B}$,
that is to the origin, we generate a partition of $\kt$ into $d_{\kt}$ triangles $\widetilde{\Delta}_i$,
$i\in\{1, 2,\ldots, d_{\kt}\}$, where $d_{\kt} \leq N_{\Tcal}$, and for which the minimum angle condition
is satisfied. The later means that there exists a constant $c_{\Tcal} > 0$, depending only on $C_{\Tcal}$ and $N_{\Tcal}$,
such that $\widetilde{h}_i(\widetilde{\rho}_i)^{-1} \leq c_{\Tcal}\quad\forall\ i\in\{1,2,\ldots,d_{\kt}\}$,
where $\widetilde{h}_i$ is the diameter of $\widetilde{\Delta}_i$ and $\widetilde{\rho}_i$ is the diameter
of the largest ball contained in $\widetilde{\Delta}_i$.
We also let $\widehat{\Delta}$ be the canonical triangle of $\R^2$ with corresponding parameters $\widehat{h}$
and $\widehat{\rho}$.
In what follows, given $K \in \Tcal_h$ and $\bzeta\in[\rmH^1(K)]^{2}$, we let $\widetilde{\bzeta}:=\bzeta\circ T_K^{-1}\in[\rmH^1(\kt)]^{2}$. With this notation  at hand, we prove the following interpolation error bound for normal components of $\rmH^1$ fuctions on edges which generalises to the VEM setting on polygons the analogous result for mixed-FEM given by Lemma 3.18 in \cite{libroGG}.
\begin{lmm}\label{lemma-approx-edges}
There exists a  constant $c_3>0$, independent of $h$, such that for all $\btau\in[\rmH^1(\Om)]^{2}$, there holds
\begin{equation}\label{approx-edges}
\Vert(\btau-\btau_I)\cdot\nu_e\Vert_{0,e}\leq c_3h_e^{1/2}|\btau|_{1,K}\quad\forall\;e\in\Ecal_h,
\end{equation}
where $K$ is any element of $\Tcal_h$ such that $K\in\omega_e$.
\end{lmm}
\begin{proof}
The proof is based on the availability of the sub-triangulation of the scaled element $\kt$ (cf. Remark~\ref{remark:meshassumptions}), and follows along the lines of the proof of Lemma 3.18 in \cite{libroGG}. Let $e\in\Ecal_h$ and $K\in\Tcal_h$ such that $K\in\omega_e$, and let $\widetilde{e}$  be the edge of $\pr\widetilde{K}$, such that $e=T_{K}^{-1}(\widetilde{e})$. We further define $T_e:=T_{K}\big\vert_e$.  Now, given $\btau\in[\rmH^1(K)]^{2}$, we know from (\ref{space-hh}) and the definition of $\btau_I$, respectively, that $\btau_I\cdot\nu\big\vert_e\in\Prm_k(e)$ and 
\begin{equation*}
\ds\int_e q(\btau-\btau_I)\cdot\nu=0\quad\forall\;q\in\Bcal_k(e).
\end{equation*}
In turns, this implies that
\begin{equation*}
\btau_I\cdot\nu_e=\Pi_k^e(\btau\cdot\nu_e),
\end{equation*}
where $\Pi_k^e:\rmL^2(e)\to\Prm_k(e)$ is the orthogonal projector. Then, it is easy to see that $\widetilde{\Pi_k^e(v)}=\Pi_{k}^{\widetilde{e}}(\widetilde{v})\quad\forall\;v\in\rmL^2(e)$, where $\Pi_{k}^{\widetilde{e}}:\rmL^2(\widetilde{e})\to\Prm_k(\widetilde{e})$ is the corresponding orthogonal projector.
Hence, we obtain
\begin{equation}\label{cta-1}
\begin{array}{lll}
\Vert(\btau-\btau_I)\cdot\nu_e\Vert_{0,e}&=&\Vert\btau\cdot\nu_e-\Pi_k^e(\btau\cdot\nu_e)\Vert_{0,e}
=\ds\frac{h_{e}^{1/2}}{h_{\widetilde{e}}^{1/2}}\Vert\widetilde{\btau\cdot{\nu}_e}-\widetilde{\Pi_k^e(\btau\cdot\nu_e)}\Vert_{0,{\widetilde{e}}}\\[3 ex]
&=&\ds\frac{h_{e}^{1/2}}{h_{\widetilde{e}}^{1/2}}\Vert\widetilde{\btau\cdot\nu}_e-\Pi_k^{\widetilde{e}}(\widetilde{\btau\cdot\nu_e})\Vert_{0,\widetilde{e}}\leq \ds\frac{h_{e}^{1/2}}{h_{\widetilde{e}}^{1/2}}\Vert\widetilde{\btau\cdot\nu_e}\Vert_{0,\widetilde{e}}\leq  \ds\frac{h_{e}^{1/2}}{h_{\widetilde{e}}^{1/2}}\Vert\widetilde{\btau}\Vert_{0,\widetilde{e}}.
\end{array}
\end{equation}
Now, let $\widetilde{\triangle}$ be the triangle formed connecting the end points of $\widetilde{e}$ to the center of $\widetilde{B}$ and  consider $\widehat{\btau}:=\widetilde{\btau}\big\vert_{\widetilde{\triangle}}\circ F\in[\rmH^1(\widehat{\triangle}]^{2}$,  where $F:\R^2\to\R^2$ is the bijective linear mapping defined by $F(\bx):=B\bx\quad\forall\;\bx\in\R^2$, with $B\in\R^{2\times 2}$ invertible, such that $F(\widehat{\triangle})$.  Let $\widehat{e}$  be the edge of $\pr\widehat{\triangle}$ such that $\widetilde{e}=F(\widehat{e})$, then
\begin{equation}\label{cta-2}
\Vert\widetilde{\btau}\Vert_{0,\widetilde{e}}=\ds\frac{h_{\widetilde{e}}^{1/2}}{h_{\widehat{e}}^{1/2}}\Vert\widehat{\btau}\Vert_{0,\widehat{e}}= \widehat{C}h_{\widetilde{e}}^{1/2}\Vert\widehat{\btau}\Vert_{0,\widehat{e}}.
\end{equation} 
Now, considering $\varphi\in C^{\infty}(\widehat{\triangle})$ such that $\varphi\equiv 1$ in a neighbourhood of $\widehat{e}$, and $\varphi\equiv 0$ in a neighbourhood of the vertex opposite to $\widehat{e}$, and applying the trace theorem in $\rmH^1(\widehat{\triangle})$, the Friedrichs-Poincar\'e inequality, and the Leibniz rule, we get
\begin{equation*}
\Vert\widehat{\btau}\Vert_{0,\widehat{e}} =\Vert\widehat{\btau}\varphi\Vert_{0,\widehat{e}}\leq \Vert\widehat{\btau}\varphi\Vert_{0,\partial\widehat{\triangle}}\leq \gamma_{\tr}\Vert\widehat{\btau}\varphi\Vert_{1,\widehat{\triangle}}\leq\gamma_{\tr} C_{p}|\widehat{\btau}\varphi|_{1,\widehat{\triangle}}\leq C_\varphi\gamma_{\tr} C_{p}|\widehat{\btau}|_{1,\widehat{\triangle}}.
\end{equation*}
Using this to bound \eqref{cta-2} and replacing the resulting bound in \eqref{cta-1} we deduce that
\begin{equation*}
\begin{array}{lll}
\Vert(\btau-\btau_I)\cdot\nu_e\Vert_{0,e} &\leq & C_\varphi\gamma_{\tr} C_{p} \widehat{C}h_{e}^{1/2}|\widehat{\btau}|_{1,\widehat{\triangle}}\leq  \widehat{C}_1C_\varphi\gamma_{\tr} C_{p} \widehat{C} h_{e}^{1/2}|\widetilde{\btau}|_{1,\widetilde{\triangle}}\\[2 ex]
&\leq & \widehat{C}_1C_\varphi\gamma_{\tr} C_{p} \widehat{C}h_{e}^{1/2}|\widetilde{\btau}|_{1,\widetilde{K}}\leq c_3 h_{e}^{1/2}|\btau|_{1,K},
\end{array}
\end{equation*}
where $c_3:=C_1\widehat{C}_1C_\varphi\gamma_{\tr} C_{p}\widehat{C} $, with $\widehat{C}_1$ and $C_1$, the $\rmH^1$-seminorm scaly constants on $\widetilde{\triangle}$ and $K$, respectively, thus concluding the proof.
\end{proof}

\subsection{A posteriori error estimator}\label{estimator}

Let $(\bsi_h,u_h)\in H_h\times Q_h$ be the unique solution of \eqref{scheme_d}.  In addition, let $\widehat{\bsi}_h,\bsi_h^{\star}$ be the discrete approximations introduced in \eqref{approx-flux} and \eqref{computo-smn-broken}, respectively. For each $K\in\Tcal_h$, we define the following local and computable error indicators:
\begin{equation*}
\begin{array}{lll}
\Phi_K^2:=\Vert f +\div\,\bsi_h\Vert_{0,K}^2,\quad
&\Lambda_{1,K}^2:=\Vert\widehat{\bsi}_h-\bsi_h^{\star}\Vert_{0,K}^2,\quad
&\Upsilon_K^2 :=\Vert(\bkappa^{-1}-\bkappa_h)\bsi_h^{\star}\Vert_{0,K}^2\,,\\[2 ex]
\Psi_K^2:=\ds\sum_{i=1}^2 \Psi_{i,K}^2\,,
&\eta_K^2:=\ds\sum_{i=1}^2\eta_{i,K}^2\,,\quad
&\theta_K^2:=\ds\sum_{i=1}^3 \theta_{i,K}^2\,,
\end{array}
\end{equation*}
where
\begin{equation*}
\begin{array}{llllll}
\Psi_{1,K}^2 &:=&\Vert\bsi_h-\widehat{\bsi}_h\Vert_{0,K}^2&\Psi_{2,K}^2 &:=&\|\bkappa^{-1}\widehat{\bsi}_h-\pin_k^0(\bkappa^{-1}\widehat{\bsi}_h)\|_{0,K}^2\\[4 ex]
\eta_{1,K}^2&:=&h_K^2\Vert\bkappa_h\bsi_h^{\star}-\nabla u_h\Vert_{0,K}^2&\eta_{2,K}^2&:=&\ds\sum_{e\in\Ecal(K)\cap\Ecal_h(\G_D)}h_e\Vert u_h-g\Vert_{0,e}^2,\\[4 ex]
\theta_{1,K}^2&:=&h_K^2\|\rot\left(\bkappa_h\bsi_h^{\star}\right)\|_{0,K}^2&\theta_{2,K}^2&:=&\ds\sum_{e\in\Ecal(K)\cap\Ecal_h(\Om)} h_e \left\Vert[\![\bkappa_h\bsi_h^{\star}\cdot s]\!]\right\Vert_{0,e}^2\,,\\[4 ex]
\theta_{3,K}^2&:=&\ds\sum_{e\in\Ecal(K)\cap\Ecal_h(\G_D)} h_e \left\Vert\bkappa_h\bsi_h^{\star}\cdot s-\frac{dg}{ds}\right\Vert_{0,e}^2\,,&&
\end{array}
\end{equation*}
and $\bkappa_h$ is a piecewise-polynomial approximation of $\bkappa^{-1}$. 

\begin{rmrk}
Notice that from the residual character  of the indicators, the computability of each local term becomes clear. 
This is the case for all terms apart from $\Psi_{1,K}$ which is not directly computable but is immediately bounded by a computable term using the stability property of Lemma~\ref{prope-ah}. 
As such, this term represents, together with $\Psi_{2,K}$, a bound on the error related to the inconsistency between the continuous and discrete bilinear forms, $a^K$ and $a_h^K$, (cf. Lemma \ref{inconsistency}).

We further observe that the last term in $\theta_{3,K}$ requires the trace $g$ to be more regular. This assumption will be stated and clarified below in Lemma \ref{cota-e2}. 
\end{rmrk}

\begin{rmrk}\label{nullterms}
If $\bkappa$ is piecewise-constant on each $K\in\Tcal_h$, we have that $\Upsilon_K$  and $\Psi_{2,K}$ are null, whereas if we use  homogeneous boundary conditions on $\G_D$, we deduce that $\eta_{2,K}$ is null.
\end{rmrk}

\begin{rmrk}
Through the a posteriori analysis below, it will be clear that the same terms but \emph{without} the postprocessing, hence with $\widehat{\bsi}_h$ in place of $\bsi_h^{\star}$ everywhere, also constitute an a posteriori bound for the error $\Vert \bsi-\bsi_h\Vert_H$. However, as we shall see, the introduction of  $\bsi_h^{\star}$ will permit us to include an optimal bound on the broken $\rmH(\div;\Om)$-norm of computable quantities.
\end{rmrk}

\subsection{Upper bound}\label{sec-rel}

We proceed with the following preliminary estimate
\begin{lmm}\label{apos-1}
Let $(\bsi,u)\in H\times Q$ and $(\bsi_h,u_h)\in H_h\times Q_h$ be the unique solutions of \eqref{scheme_c} and \eqref{scheme_d}, respectively. In addition, let $\bsi_h^{\star}$ be the discrete approximation  introduced in 
\eqref{computo-smn-broken}. Then, there exists a positive constant $C$, independent of $h$, such that
\begin{equation}\label{cota-pre}
C\Vert(\bsi,u)-(\bsi_h,u_h)\Vert_{H\times Q}\leq \left\lbrace\ds\sum_{K\in\Tcal_h}\Phi_K^2 +\Upsilon_K^2+\Psi_K^2 +\Lambda_{1,K}^2
\right\rbrace^{1/2}+\ds\sup_{\substack{\btau \in H \\ \btau{\neq}\bcero}}\frac{ E(\btau)}{\Vert\btau\Vert_{H}},
\end{equation} 
where 
\begin{equation}\label{func-e}
E(\btau):=-\ds\int_\Om\bkappa_h\bsi_h^{\star}\cdot(\btau-\btau_h)-\ds\int_\Om u_h\,\div(\btau-\btau_h)+\langle(\btau-\btau_h)\cdot\nu,g\rangle_{\G_D},
\end{equation}
for all $\btau_h\in H_h$ such that $\|\tau_h\|_{Q}\leq C\|\tau\|_H$  for some positive constant $C$ independent of $\btau$.
\end{lmm}
\begin{proof}
Consider the bounded linear operator $\Acal:H\times Q\to (H\times Q)^\prime$ induced by the  left hand-side of \eqref{scheme_c}, that is, the linear operator defined by
 \begin{equation}\label{op-apos}
 [\Acal(\brho,z),(\btau,v)]:=a(\brho,\btau)+b(\btau,z)+b(\brho,v).
\end{equation} 
From the well-posedness of the variational formulation \eqref{scheme_c},  we know that $\Acal$ is an isomorphism. In particular, there exists a positive constant $C$, such that
\begin{equation*}
C\Vert(\brho,z)\Vert_{H\times Q}\leq \ds\sup_{\substack{(\btau,v)\in {H\times Q} \\ (\btau,v){\neq}\bcero}}\frac{ [\Acal(\brho,z),(\btau,v)]}{\Vert(\btau,v)\Vert_{H\times Q}}.
\end{equation*}
Now, applying the foregoing equation to $(\brho,z):=(\bsi-\bsi_h,u-u_h)$, from \eqref{op-apos}, we get
\begin{equation}\label{esti-0}
\begin{array}{llll}
C\Vert(\bsi,u)-(\bsi_h,u_h)\Vert_{H\times Q}&\leq& \ds\sup_{\substack{(\btau,v) \in {H\times Q} \\ (\btau,v){\neq}\bcero}}\frac{ a(\bsi-\bsi_h,\btau)+b(\btau,u-u_h)+b(\bsi-\bsi_h,v)}{\Vert(\btau,v)\Vert_{H\times Q}}\\[2 ex]
&\leq& \ds\sup_{\substack{v \in Q \\ v{\neq}\bcero}}\frac{b(\bsi-\bsi_h, v)}{\Vert  v\Vert_{Q}}+\ds\sup_{\substack{\btau \in H \\ \btau{\neq}\bcero}}\frac{a(\bsi-\bsi_h,\btau)+b(\btau,u-u_h)}{\Vert\btau\Vert_{H}}\\[2 ex]
&\leq& \Vert f +\div\,\bsi_h\Vert_{0,\Om}+\ds\sup_{\substack{\btau \in H \\ \btau{\neq}\bcero}}\frac{a(\bsi-\bsi_h,\btau)+b(\btau,u-u_h)}{\Vert\btau\Vert_{H}},
\end{array}
\end{equation}
and it remains to bound the second term above.
To this end, given $\btau\in H$ and any $\btau_h\in  H_h$, from \eqref{scheme_c} and \eqref{scheme_d},  we have that
\begin{eqnarray*}
\begin{array}{llll}
& a(\bsi-\bsi_h,\btau)+b(\btau,u-u_h) =-a(\bsi_h,\btau)-b(\btau,u_h)+\langle\btau\cdot\nu,g\rangle_{\G_D}\\[2 ex]
&\quad\qquad=\langle\btau_h\cdot\nu,g\rangle_{\G_D}-a(\bsi_h,\btau)-b(\btau,u_h)+\langle(\btau-\btau_h)\cdot\nu,g\rangle_{\G_D}\\[2 ex]
&\quad\qquad=a_h(\bsi_h,\btau_h)-a(\bsi_h,\btau)-b(\btau-\btau_h,u_h)+\langle(\btau-\btau_h)\cdot\nu,g\rangle_{\G_D}\\[2 ex]
&\quad\qquad=(a_h-a)(\bsi_h,\btau_h)-a(\bsi_h-\bsi_h^{\star},\btau-\btau_h)\\[2 ex]
&\quad\qquad\quad-a(\bsi_h^{\star}-\kappa\kappa_h\bsi_h^{\star},\btau-\btau_h) \\[2 ex]
&\quad\qquad\quad-a(\kappa\kappa_h\bsi_h^{\star},\btau-\btau_h) -b(\btau-\btau_h,u_h)+\langle(\btau-\btau_h)\cdot\nu,g\rangle_{\G_D}\\[2 ex]
&\quad\qquad=:  I + II + III\,.
\end{array}
\end{eqnarray*}
Now, in what follows we take in particular $\btau_h\in H_h$  with $\|\tau_h\|_{Q}\leq C\|\tau\|_H$ for some positive constant $C$ independent of $\btau$. For $I$ and $II$, we use the bound of  Lemma \ref{inconsistency} and the Cauchy-Schwarz inequality  to deduce

\begin{equation}\label{esti-1}
I:=(a_h-a)(\bsi_h,\btau_h)-a(\bsi_h-\bsi_h^{\star},\btau-\btau_h)\leq C\left\lbrace\ds\sum_{K\in\Tcal_h}\Psi_K^2+\Lambda_{1,K}^2\right\rbrace^{1/2}\|\btau\|_H\,,
\end{equation}
and
\begin{equation}\label{esti-2}
II:=
-a(\bsi_h^{\star}-\kappa\kappa_h\bsi_h^{\star},\btau-\btau_h)\leq C\left\lbrace\ds\sum_{K\in\Tcal_h}\Upsilon_K^2\right\rbrace^{1/2}\|\btau\|_H\,,
\end{equation}

whereas bearing mind the functional $E$ (cf.\eqref{func-e}), we use the definitions \eqref{forms} to get
\begin{equation}\label{esti-3}
III:=-a(\kappa\kappa_h\bsi_h^{\star},\btau-\btau_h) -b(\btau-\btau_h,u_h)+\langle(\btau-\btau_h)\cdot\nu,g\rangle_{\G_D}= E(\btau)\,.
\end{equation}
Finally, replacing \eqref{esti-1}-\eqref{esti-3} into \eqref{esti-0}, we conclude the proof.
\end{proof}
We now aim to bound the supremum on the right hand-side of \eqref{cota-pre}, for which we need a suitable choice of $\btau_h\in H_h$ such that $\|\tau_h\|_{Q}\leq\|\tau\|_H$. To this end, in what follow we assume that the boundary $\G$ is such that $\G_N$ is contained in a convex part of $\Om$. More precisely, we make use of the following result.

\begin{lmm}\label{helm}
Assume that $\Om$ is a connected domain and that 
$\G_N$ is contained in the boundary of a convex part of $\Om$, that is there exists a convex domain $B$ such that $\Om\subset B$ and $\G_N\subseteq\pr B$. Then, for each $\btau\in H$ (cf.\eqref{space-hgn}), there exist $\bzeta\in\rmH^1(\Om)$ with $\bzeta\cdot\nu=0$ on $\G_N$ and $\chi\in\rmH^1_{\G_N}(\Om)$ (cf. Section \ref{preli}) such that
\begin{equation}\label{decom-helm}
\btau = \bzeta +\brot\,\chi\qin\Om\,,\quad\div\,\bzeta=\div\,\btau\qin\Om\,,
\end{equation}
and
\begin{equation}\label{cota-helm}
\Vert\bzeta\Vert_{1,\Om}+\Vert\chi\Vert_{1,\Om}\leq C\Vert\btau\Vert_{\div;\Om}\,,
\end{equation}
with a positive constant $C$ independent of $\btau$.
\end{lmm}
\begin{proof}
See \cite[Lemma 3.9]{mgr-2016-MMNA} for more details.
\end{proof}
Now, for $\btau\in H$ from Lemmas \ref{curl-hk} and \ref{helm}, we define $\chi_h:=\widetilde \Ical_{h}(\chi)\in \widetilde V_{h}$, and set
\begin{equation}\label{chi}
\btau_h:=\bzeta_I+\brot\,\chi_h\in H_h\,,
\end{equation}
as its associated discrete Helmholtz decomposition. Now, it follows from \eqref{chi}, the triangle inequality, \eqref{approx-pk}, \eqref{ical-1} and \eqref{cota-helm}  that 
\begin{equation*}
\|\btau_h\|_{Q}\leq \|\bzeta-\bzeta_I\|_Q+\|\bzeta\|_Q+|\chi-\chi_h|_{1,\Om}+|\chi|_{1,\Om}\leq C\|\btau\|_{H}\,,
\end{equation*}
with a positive constant $C$ independent of $\btau$. Next, we can write
\begin{equation}\label{t-chi}
\btau-\btau_h=\bzeta-\bzeta_I+\brot(\chi-\chi_h),
\end{equation}
from which, using \eqref{conmutativo}, and  the fact that $\div\,\bzeta=\div\,\btau$ in $\Om$, we deduce
\begin{equation}\label{div-t-chi}
\ds\int_\Om u_h\,\div(\btau-\btau_h)=\ds\int_\Om u_h\,\div(\bzeta-\bzeta_I)=0.
\end{equation}
Then, using the choice for $\btau_h$ given by \eqref{chi} to bound the supremum in \eqref{cota-pre}, replacing \eqref{t-chi} and \eqref{div-t-chi} into \eqref{func-e}, we find that $E(\btau)=E_1(\bzeta)+E_2(\chi)$ where
\begin{equation}\label{func-e1}
E_1(\bzeta):=-\ds\int_\Om\bkappa_h\bsi_h^{\star
}\cdot(\bzeta-\bzeta_I)+\langle(\bzeta-\bzeta_I)\cdot\nu,g\rangle_{\G_D},
\end{equation}
and
\begin{equation}\label{func-e2}
E_2(\chi):=-\ds\int_\Om\bkappa_h\bsi_h^{\star}\cdot\brot(\chi-\chi_h)+\langle\brot(\chi-\chi_h)\cdot\nu,g\rangle_{\G_D}.
\end{equation}
The following two lemmas provide the upper bounds for $|E_1(\bzeta)|$ and $|E_2(\chi)|$.
\begin{lmm}\label{cota-e1}
There exists $C>0$, independent of $h$, such that
\begin{equation*}
|E_1(\bzeta)|\leq C\left\lbrace\ds\sum_{K\in\Tcal}\eta_{K}^2\right\rbrace^{1/2}\Vert\btau\Vert_{\div;\Om}.
\end{equation*}
\end{lmm}
\begin{proof}
We  rewrite the second term in $E_1(\bzeta)$ as: 
\begin{equation}
\langle(\bzeta-\bzeta_I)\cdot\nu,g\rangle_{\G_D}=\ds\sum_{e\in\Ecal_h(\G_D)}\ds\int_e g\,(\bzeta-\bzeta_I)\cdot\nu.
\end{equation}
Next, since $u_h\big\vert_K\in\Prm_{k}(K)$, we have
\begin{equation*}
\ds\int_e u_h(\bzeta-\bzeta_I)\cdot\nu=0\quad\forall\;e\in\Ecal(K)\cap\Ecal_h(\G_D),
\end{equation*}
and
\begin{equation*}
\ds\int_K  (\bzeta-\bzeta_I)\cdot\nabla u_h=0,
\end{equation*}
for all $K\in\Tcal_h$.
Hence, using  the above expressions, we can write
\begin{equation*}
E_1(\bzeta)=-\ds\sum_{K\in\Tcal_h}\left\lbrace\ds\int_K\left(\bkappa_h\bsi_h^{\star}-\nabla u_h\right)\cdot(\bzeta-\bzeta_I)+\ds\sum_{e\in\Ecal(K)\cap\Ecal_h(\G_D)}\ds\int_e(u_h-g)(\bzeta-\bzeta_I)\cdot\nu\right\rbrace,
\end{equation*}
from which, applying the Cauchy-Schwarz inequality, the approximation properties \eqref{approx-pk} and \eqref{approx-edges}, and the fact $\Vert\bzeta\Vert_{1,\Om}\leq \Vert\btau\Vert_{\div;\Om}$, we obtain the required estimate.
\end{proof}
\begin{lmm}\label{cota-e2}
Assume that $g\in\rmH^1(\G_D)$. Then, there exists $C>0$, independent of $h$, such that
\begin{equation*}
|E_2(\chi)|\leq C\left\lbrace\ds\sum_{K\in\Tcal}\theta_K^2\right\rbrace^{1/2}\Vert\btau\Vert_{\div;\Om}.
\end{equation*}
\end{lmm}
\begin{proof}
We proceed as in the proof of the Lemma 3.11 in \cite{mgr-2016-MMNA}. Integrating by parts on each $K\in\Tcal_h$, using that $\brot(\chi-\chi_h)\cdot\nu=\ds\frac{d}{ds}(\chi-\chi_h)$, noting that $\ds\frac{dg}{ds}\in\rmL^2(\G_D)$, and  using  the  fact  that $\chi\big\vert_{\G_N}=\chi_h\big\vert_{\G_N}=0$, we get
\begin{equation*}
\begin{array}{lll}
E_2(\chi)&=&-\ds\sum_{K\in\Tcal_h}\ds\int_K\bkappa_h\bsi_h^{\star}\cdot\brot(\chi-\chi_h)+\left\langle\frac{d}{ds}(\chi-\chi_h),g\right\rangle_{\G_D}\\[2 ex]
&=&-\ds\sum_{K\in\Tcal_h}\Big\{\ds\int_{K}\rot\left(\bkappa_h\bsi_h^{\star}\right)(\chi-\chi_h)-\ds\int_{\pr K}\left(\bkappa_h\bsi_h^{\star}\cdot s_K\right)(\chi-\chi_h)\Big\} -\ds\int_{\G_D}\frac{dg}{ds}(\chi-\chi_h)\\[2 ex]
 &=&-\ds\sum_{K\in\Tcal_h}\Big\{\ds\int_{K}\rot\left(\bkappa_h\bsi_h^{\star}\right)(\chi-\chi_h)-\ds\sum_{e\in\Ecal(K)\cap\Ecal_h(\Om)}\ds\int_{e}[\![\bkappa_h\bsi_h^{\star}\cdot s]\!](\chi-\chi_h) \quad\\[2 ex]
& &\qquad\qquad -\ds\sum_{e\in\Ecal(K)\cap\Ecal_h(\G_D)}\ds\int_{e}\left(\bkappa_h\bsi_h^{\star}\cdot s-\frac{dg}{ds}\right)(\chi-\chi_h)\Big\}.
\end{array}
\end{equation*}
In this way, since $\chi_h=\widetilde \Ical_{h}(\chi)$, applying the Cauchy-Schwarz inequality to each term in the above expression and making use of the approximation properties~\eqref{ical-1} and~\eqref{ical-2} and the fact that the number of elements in $\omega_e$ is bounded, we conclude the proof.
\end{proof}

Finally, from Lemmas \ref{apos-1}, \ref{cota-e1} and \ref{cota-e2} we deduce an upper bound for the global error.
\begin{thrm}\label{main0-thm}
Let $(\bsi,u)\in H\times Q$ and $(\bsi_h,u_h)\in H_h\times Q_h$ be the unique solutions of the problem \eqref{scheme_c} and \eqref{scheme_d}, respectively. Then, there exists a positive constant $C$, independent of $h$, such that 
\begin{equation*}
\Vert(\bsi,u)-(\bsi_h,u_h)\Vert_{H\times Q}\leq C
\left\lbrace\ds\sum_{K\in\Tcal_h} 
\Phi_K^2
+\Upsilon_K^2+\Psi_K^2
+\Lambda_{1,K}^2
+\eta_K^2+\theta_K^2 
\right\rbrace^{1/2}.
\end{equation*}
\end{thrm}

We recall from the discussion in Section~\ref{sec:computable} that the corresponding result for the computable quantity $\widehat{\bsi}_h$ is only to be expected in the $\rmL^2$-norm. Instead, 
 for the error using the postprocessing flux we are able to obtain the following result in line with Theorem~\ref{main0-thm}.

\begin{thrm}\label{main1-thm}
Let $(\bsi,u)\in H\times Q$ and $(\bsi_h,u_h)\in H_h\times Q_h$ be the unique solutions of the problem \eqref{scheme_c} and \eqref{scheme_d}, respectively. In addition, let $\bsi_h^{\star}$ be the discrete postprocessing  introduced in 
\eqref{computo-smn-broken}. Then, there exists a positive constant $C$, independent of $h$, such that 
\begin{equation*}
\left\lbrace\sum_{K\in\Tcal_h}\|\bsi-\bsi_{h,K}^{\star}\|_{\div;K}^2\right\rbrace^{1/2} +\| u-u_h\|_{Q}\leq C
\left\lbrace\ds\sum_{K\in\Tcal_h} 
\Phi_K^2
+\Upsilon_K^2+\Psi_K^2
+\Lambda_{K}^2
+\eta_K^2+\theta_K^2 
\right\rbrace^{1/2},
\end{equation*}
with
\[
\Lambda_K^2:=\ds\sum_{i=1}^2\Lambda_{i,K}^2\quad\text{where}\quad
\Lambda_{2,K}^2:=\Vert\div\,\bsi_h-\div\,\bsi_h^{\star}\Vert_{0,K}^2.
\]
\end{thrm}
\begin{proof}
From the triangle inequality, we have
\begin{equation*}
\begin{array}{lll}
\|\bsi-\bsi_{h,K}^{\star}\|_{\bdiv;K} &\leq & \|\bsi-\bsi_h\|_{\div;K} + \|\bsi_h-\bsi_{h,K}^{\star}\|_{0,K}+\|\div\,\bsi_h-\div\,\bsi_{h,K}^{\star}\|_{0,K}\\[2 ex]
&\leq & \|\bsi-\bsi_h\|_{\div;K} + \|\bsi_h-\widehat{\bsi}_h\|_{0,K}+\|\widehat{\bsi}_h-\bsi_{h,K}^{\star}\|_{0,K}+\|\div\,\bsi_h-\div\,\bsi_{h,K}^{\star}\|_{0,K}\,.
\end{array}
\end{equation*}
Then, since $\rmH(\div;\Om)\subset\rmH(\div;\Tcal_h)$ and using the definition of $\Psi_K^2$ and $\Lambda_K^2$, we get
\begin{equation*}
\left\lbrace\sum_{K\in\Tcal_h}\|\bsi-\bsi_{h,K}^{\star}\|_{\div;K}^2\right\rbrace^{1/2}\leq C\left\lbrace\|\bsi-\bsi_h\|_H+\left\lbrace\ds\sum_{K\in\Tcal_h}\Psi_K^2+\Lambda_K^2\right\rbrace^{1/2}\right\rbrace\,.
\end{equation*}
Threrefore, the result is consequence of the foregoing equation and the Theorem \ref{main0-thm}.
\end{proof}

\subsection{Lower bound}\label{sec-eff}

In this section we derive suitable upper bounds for the terms defining the local error indicators. First, using that $f=-\div\,\bsi$ in $\Om$ we have that
\begin{equation}\label{cota-phik}
\Phi_K^2= \Vert\div(\bsi-\bsi_h)\Vert_{0,K}^2\leq 2\Big\{\|\bsi-\bsi_h^{\star}\|_{\div;K}^2+\Lambda_{1,K}^2\Big\}\,.
\end{equation}
Moreover, adding and subtracting $\bsi$, we easily have

\begin{equation}\label{cota-lamdak}
\begin{array}{ll}
\Lambda_K^2&=\Vert\widehat{\bsi}_h-\bsi_h^{\star}\Vert_{0,K}^2+\Vert\div\,\bsi_h-\div\,\bsi_h^{\star}\Vert_{0,K}^2\\[2 ex]
&\leq 2\left\{
\Vert\bsi-\widehat{\bsi}_h\Vert_{0,K}^2+\|\bsi-\bsi_h^{\star}\|_{\div;K}^2+\Phi_K^2
\right\}\,.
\end{array}
\end{equation}

In addition, proceeding as in \cite[Lemma 18]{cgto-2017-NM}, we deduce 
\begin{equation}\label{cota-eta3k}
\Psi_{2,K}^2\leq C\Big\{\Lambda_{1,K}^2+\|\bsi-\bsi_h^{\star}\|_{0,K}^2 + \|\bkappa^{-1}\bsi-\pin_k^0(\bkappa^{-1}\bsi)\|_{0,K}^2\Big\}\,,
\end{equation}
with $C$ depending only on $\bkappa$ and $\widehat{c}_0$. 

\begin{rmrk}
Again by adding and subtracting $\bsi$ we have
\begin{equation}\label{cota-eta4k}
\Psi_{1,K}^2=\Vert\bsi_h-\widehat{\bsi}_h\Vert_{0,K}^2
\leq 2\left\{\Vert\bsi-\bsi_h\Vert_{0,K}^2+\Vert\bsi-\widehat{\bsi}_h\Vert_{0,K}^2 \right\}.
\end{equation}
This does provide a lower bound, although in terms of the error $\bsi-\bsi_h$. Here we have chosen, instead,  to leave this  term as is, interpreting it as a sort of oscillation term representing the virtual inconsistency of the method.
\end{rmrk}

The upper bounds of the terms which depend on the mesh parameters $h_K$ and $h_e$, will be derived next. To this end, we proceed similarly as in \cite{c-1997-MC} and \cite{c-1998-NM} and apply the technique based on bubble functions, together with inverse inequalities. Following \cite[Section 4]{cgto-2017-NM} and \cite[Section 3]{mrr-2017-CMA}, given $K\in\Tcal_h$, a  bubble  function $\psi_K$ can be constructed piecewise as the sum of the (polynomial) barycentric bubble functions (cf.\cite{verfurth,mt-2000-PAM}) on each triangle of the shape-regular sub-triangulation of the
mesh element $K$ discussed in the Section \ref{proper}. Further, an edge bubble function $\psi_e$, $e\in\pr K$, is a piecewise quadratic function attaining the value 1 at the mid-point of $e$ and vanishing on the triangles that do not contain $e$ on their boundary. Furthermore, given $k\geq 0$, there exists an extension operator $L:\rmC(e)\to\rmC(K)$ that satisfies $L(p)\in\Prm_k(K)$ and $L(p)\big\vert_e=p$ for all $p\
 \in\Prm_k(e)$ (cf.\cite[Remark 3.1]{mrr-2017-CMA}). Further properties of $\psi_K,\psi_e$, and $L$ are stated in  the following lemma. See \cite[Section 4]{cgto-2017-NM} and \cite[Section 3]{mrr-2017-CMA} for more details.
\begin{lmm}\label{cota-burlbles}
Given $k\geq 0$ and $K\in\Tcal_h$, there exists a positive constant $C_{\mathrm{bub}}$, independent of $h_K$ such that 
\begin{equation}\label{burble-1}
C_{\mathrm{bub}}^{-1}\Vert q\Vert_{0,K}^2\leq \Vert\psi_K^{1/2}q\Vert_{0,K}^2\leq C_{\mathrm{bub}}\Vert q\Vert_{0,K}^2\qquad\forall\;q\in\Prm_k(K),
\end{equation}
and
\begin{equation}\label{burble-2}
C_{\mathrm{bub}}^{-1}\Vert q\Vert_{0,K}\leq\Vert\psi_K q\Vert_{0,K}+h_K|\psi_K q|_{1,K}\leq C_{\mathrm{bub}}\Vert q\Vert_{0,K}\qquad\forall\;q\in\Prm_k(K).
\end{equation}
In addition, given $e\in\pr K$, there hold
\begin{equation}\label{burble-3}
C_{\mathrm{bub}}^{-1}\Vert q\Vert_{0,e}^2\leq \Vert\psi_e^{1/2}q\Vert_{0,e}^2\leq C_{\mathrm{bub}}\Vert q\Vert_{0,e}^2\qquad\forall\;q\in\Prm_k(e),
\end{equation}
and
\begin{equation}\label{burble-4}
h_K^{-1/2}\Vert \psi_e L(q)\Vert_{0,K}+ h_K^{1/2}|\psi_e L(q)|_{1,K}\leq C_{\mathrm{bub}}\Vert q\Vert_{0,e}\qquad\forall\;q\in\Prm_k(e),
\end{equation}
where $K\in\omega_e$.
\end{lmm}

We start the analysis bounding the terms defining $\eta_{1,K}^2$ and $\eta_{2,K}^2$.
\begin{lmm}\label{cota-pcaln-nuh}
There exists a constant $C>0$, independent of $h$, such that
\begin{equation*}
h_K^2\|\bkappa_h\bsi_h^{\star}-\nabla u_h\|_{0,K}^2\leq C\Big\{h_K^{2}\| \bsi-\bsi_h^{\star}\|_{0,K}^2+ h_K^2\|(\bkappa^{-1}-\bkappa_h)\bsi_h^{\star}\|_{0,K}^2+\| u-u_h\|_{0,K}^2\Big\}\qquad\forall\; K\in\Tcal_h.
\end{equation*}
\end{lmm}
\begin{proof}
It is a slight modification of the proof of the Lemma 6.3 in \cite{c-1997-MC} (see also Lemma 5.5 in \cite{g-2004-ETNA}). Given $K\in\Tcal_h$ we denote $\gamma_K:=\bkappa_h\bsi_h^{\star}-\nabla u_h\in[\Prm_\ell(K)]^{2}$ for some $\ell\geq 0$. Then, applying \eqref{burble-1}, using that $\bkappa^{-1}\bsi=\nabla u$ in $\Om$, and integrating by parts, we find that
\begin{equation*}
\begin{array}{lll}
 C_{\mathrm{bub}}^{-1}\Vert\gamma_K\Vert_{0,K}^2&\leq& \Vert\psi_K^{1/2}\gamma_K\Vert_{0,K}^2=\ds\int_K\psi_K\gamma_K\cdot\left\lbrace\bkappa_h\bsi_h^{\star}-\nabla u_h\right\rbrace\\[2 ex]
&=&\ds\int_K\psi_K\gamma_K\cdot\left\lbrace(\bkappa_h-\bkappa^{-1})\bsi_h^{\star}+\bkappa^{-1}\bsi_h^{\star}-\bkappa^{-1}\bsi+\nabla(u-u_h)\right\rbrace\\[2 ex]
&=&-\ds\int_K\psi_K\gamma_K\cdot(\bkappa^{-1}-\bkappa_h)\bsi_h^{\star}
-\ds\int_K\psi_K\gamma_K\cdot\bkappa^{-1}(\bsi-\bsi_h^{\star})\\[2 ex]
&&-\ds\int_K\div(\psi_K\gamma_K)(u-u_h).
\end{array}
\end{equation*}
Then, applying the Cauchy-Schwarz inequality, the estimate \eqref{burble-2}, and setting $C_\bkappa:=\max\Big\{1, \|\bkappa^{-1}\|\Big\}$, we get
\begin{equation*}
\begin{array}{lll}
 C_{\mathrm{bub}}^{-1}\|\gamma_K\|_{0,K}^2 &\leq &  C_\bkappa\Big\{\|\psi_K\gamma_K\|_{0,K}\Big\{\|(\bkappa^{-1}-\bkappa_h)\bsi_h^{\star}\|_{0,K}+\|\bsi-\bsi_h^{\star}\|_{0,K}\Big\} +|\psi_K\gamma_K|_{1,K}\| u-u_h\|_{0,K}\Big\}\\[2 ex]
&\leq&C_\bkappa C_{\mathrm{bub}}\left\lbrace\|(\bkappa^{-1}-\bkappa_h)\bsi_h^{\star}\|_{0,K}+\|\bsi-\bsi_h^{\star}\|_{0,K}+h_K^{-1}\| u-u_h\|_{0,K}\right\rbrace\|\gamma_K\|_{0,K}\\[2 ex]
&\leq & 2C_\bkappa C_{\mathrm{bub}}\left\lbrace\|(\bkappa^{-1}-\bkappa_h)\bsi_h^{\star}\|_{0,K}^2+\|\bsi-\bsi_h^{\star}\|_{0,K}^2+h_K^{-2}\| u-u_h\|_{0,K}^2\right\rbrace^{1/2}\|\gamma_K\|_{0,K},
\end{array}
\end{equation*}
whence, the proof is concluded.
\end{proof}
\begin{lmm}
There exists a constant $C>0$, independent of $h$, such that
\begin{equation*}
h_e\| g-u_h\|_{0,e}^2\leq C\Big\{h_K^2\|\bsi-\bsi_h^{\star}\|_{0,K}^2+ h_K^2\|(\bkappa^{-1}-\bkappa_h)\bsi_h^{\star}\|_{0,K}^2+\Vert u-u_h\Vert_{0,K}^2\Big\}\quad\forall\, e\in\Ecal_h(\G_D),
\end{equation*}
where $K\in\Tcal_h$ is such that $e\in\pr K$.
\end{lmm}
\begin{proof}
We proceed as in the proof of the Lemma 4.14 in \cite{gms-2010-CMAME}. We consider $e\in\Ecal_h(\G_D)$ and $K\in\Tcal_h$ such that $e\in\pr K$. Then, applying a trace inequality, together with the fact that $u=g$ on $\G_D$ and $\bsi=\bkappa\nabla u$ in $\Om$, we get
\begin{equation*}
\begin{array}{lll}
& \| g-u_h\|_{0,e}^2=\| u-u_h\|_{0,e}^2\leq C_\tr\left\lbrace h_K^{-1}\| u-u_h\|_{0,K}^2+h_K|u-u_h|_{1,K}^2\right\rbrace\\[2ex]
&\leq  2C_\bkappa C_\tr\Big\{ h_K^{-1}\| u-u_h\|_{0,K}^2+h_K\Big\{\| \bsi-\bsi_h^{\star}\|_{0,K}^2+\|\bkappa_h\bsi_h^{\star}-\nabla u_h\|_{0,K}^2+\|(\bkappa^{-1}-\bkappa_h)\bsi_h^{\star}\|_{0,K}^2\Big\}\Big\},
\end{array}
\end{equation*}
with $C_\bkappa $ as in the proof of Lemma~\ref{cota-pcaln-nuh}.

From this, using the bound $h_e\leq h_K$ and the estimate of Lemma \ref{cota-pcaln-nuh} we obtain the result.
\end{proof}

The following result is required in view of proving upper bounds for the terms defining $\theta_K^2$.
\begin{lmm}\label{cotas-rot-jump}
Let $\bzeta_h\in[\rmL^2(\Om)]^{2}$ be a piecewise polynomial of degree $k\geq 0$ on each $K\in\Tcal_h$. In addition let $\bzeta\in[\rmL^2(\Om)]^{2}$ be such that $\rot(\bzeta)=0$ in $\Om$. Then, there exists $C>0$, depending only on $C_{\mathrm{bub}}$, such that
\begin{equation}\label{cota-bucurl}
\|\rot\,\bzeta_h\|_{0,K}\leq Ch_K^{-1}\|\bzeta-\bzeta_h\|_{0,K}\quad\forall\,K\in\Tcal_h\,,
\end{equation}
and
\begin{equation}\label{cota-jump}
\Vert[\![\bzeta_h\cdot s_e]\!]\Vert_{0,e}\leq Ch_e^{-1/2}\Vert\bzeta-\bzeta_h\Vert_{0,\omega_e}\quad\forall\;e\in\Ecal_h(\Om)\,.
\end{equation}
\end{lmm}
\begin{proof}
To show~\eqref{cota-bucurl}, we proceed as in the proof of Lemma 4.3 in  \cite{bgg-2006-M2AN}. Applying \eqref{burble-1}, observing that $\psi_K=0$ on $\pr K$, and using the Cauchy-Schwarz inequality, we get
\begin{equation*}
\begin{array}{lll}
C_{\mathrm{bub}}^{-1}\|\rot\,\bzeta_h\|_{0,K}^2 &\leq & \|\psi_K^{1/2}\rot\,\bzeta_h\|_{0,K}^2=-\ds\int_K\psi_K\rot\,\bzeta_h\,\rot(\bzeta-\bzeta_h)\\[2 ex]
 & =& \ds\int_K(\bzeta-\bzeta_h)\cdot\brot(\psi_K\rot\,\bzeta_h)\leq \|\bzeta-\bzeta_h\|_{0,K}|\psi_K\rot\,\bzeta_h|_{1,K}\,.
\end{array}
\end{equation*}
Then, from inverse inequality \eqref{burble-2}, we deduce \eqref{cota-bucurl}.

The estimate~\eqref{cota-jump} follows from a slight
modification of the proof of \cite[Lemma 4.4]{bgg-2006-M2AN}. Indeed,  given $e\in\Ecal_h(\Om)$, we let $J_h:=[\![\bzeta_h\cdot s_e]\!]\in\Prm_k(e)$. Then, utilizing \eqref{burble-3}, the fact that $[\![\bzeta\cdot s_e]\!]=0$ a.e on $e$, and integrating by parts on each $K\in\Tcal_h$, we get
\begin{equation*}
\begin{array}{lll}
C_{\mathrm{bub}}^{-1}\Vert J_h\Vert_{0,e}^2 &\leq & \Vert\psi_e^{1/2}J_h\Vert_{0,e}^2=\Vert\psi_e^{1/2}L(J_h)\Vert_{0,e}^2=\ds\int_{e}\psi_e L(J_h)[\![\bzeta_h\cdot s]\!]\\[2 ex]
&=&\ds\int_{\omega_e}(\bzeta_h-\bzeta)\cdot\brot(\psi_e L(J_h)) + \ds\int_{\omega_e} \psi_e L(J_h)\,\rot\,\bzeta_h\,,
\end{array}
\end{equation*}
which, using the Cauchy-Schwarz inequality, the estimates \eqref{burble-4} and \eqref{cota-bucurl}, and the fact that $h_e\leq h_K$, yields
\begin{equation*}
\begin{array}{lll}
C_{\mathrm{bub}}^{-1}\Vert J_h\Vert_{0,e}^2 &\leq & |\psi_e L(J_h)|_{1,\omega_e}\Vert\bzeta-\bzeta_h\Vert_{0,\omega_e}+\|\psi_e L(J_h)\|_{0,\omega_e}\|\rot\,\bzeta_h\|_{0,\omega_e}\\[2 ex]
&\leq & 2N_\Tcal C_{\mathrm{bub}}h_e^{-1/2}\Vert\bzeta-\bzeta_h\Vert_{0,\omega_e}\Vert J_h\Vert_{0,e}\,,
\end{array}
\end{equation*}
whence, we conclude the proof of \eqref{cota-jump}.
\end{proof}

\begin{lmm}\label{cotas-rot-pk}
There exists $C>0$, independent of $h$, such that
\begin{equation*}
h_K^2\|\rot(\bkappa_h\bsi_h^{\star})\|_{0,K}^2\leq C\Big\{\|\bsi-\bsi_h^{\star}\|_{0,K}^2 + \|(\bkappa^{-1}-\bkappa_h)\bsi_h^{\star}\|_{0,K}^2\Big\}\quad\forall\, e\in\Ecal_h(\Om),
\end{equation*}
and
\begin{equation*}
h_e\Vert[\![\bkappa_h\bsi_h^{\star}\cdot s_e]\!]\Vert_{0,e}^2\leq C\Big\{\|\bsi-\bsi_h^{\star}\|_{0,K}^2 + \|(\bkappa^{-1}-\bkappa_h)\bsi_h^{\star}\|_{0,K}^2\Big\}\quad\forall\, e\in\Ecal_h(\Om),
\end{equation*}
where $K\in\Tcal_h$ is such that $K\in\omega_e$.
\end{lmm}
\begin{proof}
It suffices to apply Lemma \ref{cotas-rot-jump} with $\bzeta_h:=\bkappa_h\bsi_h^{\star}$ and $\bzeta:=\bkappa^{-1}\bsi=\nabla u$, and the triangle inequality.
\end{proof}
\begin{lmm}\label{cota-g}
Assume that $\frac{dg}{ds}$ is piecewise polynomial on $\G_D$.
Then, there exists $C>0$, independent of $h$, such that 

\begin{equation}\label{esti-cota-g}
h_e \left\Vert\bkappa_h\bsi_h^{\star}\cdot s-\frac{dg}{ds}\right\Vert_{0,e}^2\leq C\Big\{\|\bsi-\bsi_h^{\star}\|_{0,K}^2+\|(\bkappa^{-1}-\bkappa_h)\bsi_h^{\star}\|_{0,K}^2\Big\}\quad\forall\,e\in\Ecal_h(\G_D),
\end{equation}

where $K\in\Tcal_h$ is such that $K\in\omega_e$.
\end{lmm}
\begin{proof}
We proceed as in the proof of  Lemma 4.15 in \cite{gms-2010-CMAME} (see also Lemma 5.7 in \cite{g-2004-ETNA}). Given $e\in\Ecal_h(\G_D)$  and $K\in\omega_e$, we denote $\gamma_e:=\bkappa_h\bsi_h^{\star}\cdot s-\ds\frac{dg}{ds}\in \Prm_\ell(e)$ for some $\ell\geq 0$. Then, applying \eqref{burble-3}, the fact that $\nabla u\cdot s=\ds\frac{dg}{ds}$, integrating by parts and using that $\bkappa^{-1}\bsi=\nabla u$ in $\Om$, we obtain that
\begin{equation*}
\begin{array}{lll}
 C_{\mathrm{bub}}^{-1}\Vert\gamma_e\Vert_{0,e}^2 &\leq & \Vert\psi_e^{1/2}\gamma_e\Vert_{0,e}^2=\ds\int_e\psi_e\gamma_e\Big\{\bkappa_h\bsi_h^{\star}\cdot s-\nabla u\cdot s\Big\}\\[2 ex]
&=&-\ds\int_{\pr K}\psi_{e} L(\gamma_e)\left\lbrace\left(\bkappa^{-1}\bsi-\bkappa_h\bsi_h^{\star}\right)\cdot s\right\rbrace\\[2 ex]
&=&-\ds\int_{K}\brot(\psi_{e} L(\gamma_e))\cdot(\bkappa^{-1}-\bkappa_h)\bsi_h^{\star} -\ds\int_{K}\brot(\psi_{e} L(\gamma_e))\cdot\bkappa^{-1}\Big\{\bsi-\bsi_h^{\star}\Big\} \\[2 ex]
&& + \ds\int_K\psi_e L(\gamma_e)\rot(\bkappa_h\bsi_h^{\star})\,.
\end{array}
\end{equation*}
Next, applying the Cauchy-Schwarz inequality, Lemma \ref{cotas-rot-pk}, the estimate \eqref{burble-4}, and the fact that $h_e\leq h_K$ we get
\begin{equation*}
\begin{array}{lll}
 C_{\mathrm{bub}}^{-1}\|\gamma_e\|_{0,e}^2&\leq& 
C\Big\{|\psi_{e} L(\gamma_e)|_{1,K}+h_K^{-1}\|\psi_{e} L(\gamma_e)\|_{0,K}\Big\}\Big\{\|\bsi-\bsi_h^{\star}\|_{0,K}+\|(\bkappa^{-1}-\bkappa_h)\bsi_h^{\star}\|_{0,K}\Big\}\\[2 ex]
&\leq& Ch_e^{-1/2}\Big\{\|\bsi-\bsi_h^{\star}\|_{0,K}+\|(\bkappa^{-1}-\bkappa_h)\bsi_h^{\star}\|_{0,K}\Big\}\|\gamma_e\|_{0,e},
\end{array}
\end{equation*}
and the proof is complete.
\end{proof}

If $\frac{dg}{ds}$  is not piecewise polynomial but sufficiently smooth, Lemma~\ref{cota-g} can still be proven with higher order terms given by the errors arising from suitable polynomial approximations appearing in \eqref{esti-cota-g}.

Finally, a lower bound is obtained from estimates \eqref{cota-phik}-\eqref{cota-eta3k}, together with Lemmata~\ref{cota-pcaln-nuh} throughout \ref{cota-g}, after summing up over $K\in\Tcal_h$ and using the fact that the number of elements on each domain $\omega_e$ is bounded.

\section{Numerical Tests}\label{sec:numerics}

In this section, we present three numerical tests confirming the upper and  lower bounds, derived in Section~\ref{sec:aposteriorianalysis}, for the a posteriori
error estimator of Theorem~\ref{main1-thm}, and showing the behaviour of the associated adaptive algorithm.  We begin by introducing additional notations. In what follows, $N$ stands for the total number of degrees of freedom of \eqref{scheme_d}, that is,
\begin{equation*}
N:=(k+1)\times\lbrace\mbox{number of edges }e\in\Tcal_h\rbrace+\ds\frac{(k+2)(3k+1)}{2}\times
\lbrace\mbox{number of elements }K\in\Tcal_h\rbrace.
\end{equation*}
Also, the individual errors are defined by
\begin{equation*}
{\tt e}(\bsi)\, :=\, \left\lbrace\sum_{K\in\Tcal_h}\|\bsi-\bsi_{h,K}^{\star}\|_{\div;K}^2\right\rbrace^{1/2},\quad {\tt e}(u)\, :=\, \|u-u_h\|_{0,\Om},\qy{\tt e}(\bsi,u)\, :=\,\Big\{ [{\tt e}(\bsi)]^2 + [{\tt e}(u)]^2\Big\}^{1/2}\,,
\end{equation*}
whereas
the associated experimental rates of convergence are given by
\begin{equation*}
{\tt r}(\cdot):=-2\frac{\log({\tt e}(\cdot)/{\tt e}^\prime(\cdot))}{\log(N/N^\prime)},
\end{equation*}
where ${\tt e}$ and ${\tt e}^\prime$ denote the corresponding errors for two consecutive meshes with  $N$ and $N^\prime$ denote the corresponding degrees of freedom of each decomposition. Denote by $\Theta$ the a posteriori error estimator of Theorem~\ref{main1-thm}. The effectivity  of the estimator $\Theta$ is given by
\begin{equation*}
{\tt eff}(\Theta):=\dfrac{{\tt e}(\bsi,u)}{\Theta}\,.
\end{equation*}
For the tests that include adaptivity, we use the strategy:
\begin{enumerate}
\item [(i)] Start with a coarse mesh $\Tcal_h$.
\item [(ii)] Solve the discrete problem on the current mesh $\Tcal_h$.
\item [(iii)] Compute local indicators for each $K\in\Tcal_h$.
\item  [(iv)] Mark each $K^\prime\in\Tcal_h$ such that
\begin{equation*}
\Theta_{K^\prime}\geq\beta\ds\max_{K\in\Tcal_h}\Theta_K,
\end{equation*}
  with $\beta\in [0,1]$ and we refine using the midpoint of each edge of each element and connecting this to its barycentre. Here, we use $\beta = 0.5$.
\item [(v)]  Update $\Tcal_h$ with the new mesh and go to step (ii).
\end{enumerate}
Hereafter, in all numerical tests we have $\bkappa =\begin{pmatrix}
1 & 0\\
0 & 1
\end{pmatrix}$ and we consider domains $\Om$  satisfying Lemma~\ref{helm}. In this case, we have that $\Upsilon_K^2$ and $\Psi_{2,K}^2$ are  null for each $K\in\Tcal_h$ (cf.  Remark~\ref{nullterms} in Section~\ref{estimator}).

\subsection{Test 1. Smooth solution: behaviour of the estimator under uniform refinement}

For this test case, we consider $\Omega:=(0,1)^2$ with $\G_D:=\Big\{(w,0),(0,w)\in\Omega:\quad 0\leq w\leq 1 \Big\}$ and $\G_N:=\G\setminus\overline{\G}_D$. The source term $f$ and
the boundary data $g$ are chosen such that the exact solution is given by $u(x,y) =\cos(\pi x)\cos(\pi y)$

\begin{table}
\centering
\scalebox{0.88}{\begin{tabular}{|c|r|c c|c c|c c|c c|c|}\hline
$k$& $N$ & ${\tt e}(\bsi)$ & ${\tt r}(\bsi)$ & ${\tt e}(u)$ & ${\tt r}(u)$  & ${\tt e}(\bsi,u)$ & ${\tt r}(\bsi,u)$  & $\Theta$ & ${\tt r}(\Theta)$ & ${\tt eff}(\Theta)$\\ \hline
  &    589 &  1.0959e+00 &  $--$  &  5.5406e-02 &   $--$  & 1.0973e+00 &  $--$  & 1.2946e+00 &   $--$  &  8.4760e-01 \\
  &   3469 &  4.3834e-01 & 1.0336 &  2.1963e-02 &  1.0437 & 4.3889e-01 & 1.0336 & 5.1997e-01 &  1.0288 &  8.4405e-01 \\
0 &   8749 &  2.7388e-01 & 1.0167 &  1.3708e-02 &  1.0191 & 2.7423e-01 & 1.0167 & 3.2511e-01 &  1.0153 &  8.4348e-01 \\
  &  19605 &  1.8247e-01 & 1.0067 &  9.1290e-03 &  1.0077 & 1.8270e-01 & 1.0067 & 2.1667e-01 &  1.0059 &  8.4323e-01 \\
  &  43805 &  1.2165e-01 & 1.0087 &  6.0851e-03 &  1.0090 & 1.2180e-01 & 1.0087 & 1.4447e-01 &  1.0082 &  8.4307e-01 \\ \hline
  &   1766 &  5.9382e-02 &  $--$  &  3.0274e-03 &   $--$  & 5.9459e-02 &  $--$  & 7.2082e-02 &   $--$  &  8.2488e-01 \\
  &  10406 &  9.5820e-03 & 2.0569 &  4.8269e-04 &  2.0704 & 9.5941e-03 & 2.0569 & 1.1588e-02 &  2.0611 &  8.2791e-01 \\
1 &  26246 &  3.7498e-03 & 2.0282 &  1.8859e-04 &  2.0317 & 3.7546e-03 & 2.0282 & 4.5310e-03 &  2.0301 &  8.2863e-01 \\
  &  58814 &  1.6674e-03 & 2.0088 &  8.3810e-05 &  2.0103 & 1.6695e-03 & 2.0088 & 2.0135e-03 &  2.0105 &  8.2918e-01 \\
  & 131414 &  7.4166e-04 & 2.0154 &  3.7268e-05 &  2.0160 & 7.4260e-04 & 2.0154 & 8.9537e-04 &  2.0159 &  8.2938e-01 \\ \hline
  &   3384 &  2.2263e-03 &  $--$  &  1.1132e-04 &   $--$  & 2.2290e-03 &  $--$  & 3.6826e-03 &   $--$  &  6.0529e-01 \\
  &  19944 &  1.4599e-04 & 3.0718 &  7.2751e-06 &  3.0757 & 1.4618e-04 & 3.0718 & 2.3901e-04 &  3.0835 &  6.1158e-01 \\
2 &  50304 &  3.5848e-05 & 3.0358 &  1.7855e-06 &  3.0368 & 3.5892e-05 & 3.0358 & 5.8541e-05 &  3.0412 &  6.1312e-01 \\
  & 112726 &  1.0645e-05 & 3.0097 &  5.3011e-07 &  3.0101 & 1.0658e-05 & 3.0097 & 1.7356e-05 &  3.0135 &  6.1405e-01 \\
  & 251876 &  3.1615e-06 & 3.0200 &  1.5743e-07 &  3.0202 & 3.1654e-06 & 3.0200 & 5.1503e-06 &  3.0222 &  6.1460e-01 \\ \hline  
\end{tabular} }\vspace{0.5cm}   
\caption{Test 1. Convergence history for an uniformly generated sequence of  hexagonal meshes.}  
\label{tab1}  
\end{table}

Table \ref{tab1} shows the convergence history of the error for each variable and the estimator on a sequence of uniformly refined hexagonal meshes, indicating that both converge at the optimal rate for polynomial degrees $k = 0,1,2$. Moreover, the effectivity remains bounded. In addition, we see from Table \ref{tab2} that each term of the error estimator converge with optimal order $k+1$.

\begin{table}
\centering
\scalebox{0.84}{\begin{tabular}{|c|r|c c|c c|c c|c c|c c|c c|c c|c c|}\hline
$k$ &$N$ & $\Phi$ & ${\tt r}(\Phi)$ & $\eta$ & ${\tt e}(\eta)$ & $\theta$  & ${\tt r}(\theta)$ & $\Psi$  & ${\tt r}(\Psi)$ & $\Lambda$ & ${\tt r}(\Lambda)$ \\ \hline
  &   589 & 1.0813e+00 &  $--$  & 2.4480e-01 &  $--$  &  4.1938e-01 &  $--$  &  4.9543e-01 &  $--$  & 1.5941e-01  &  $--$  \\
  &  3469 & 4.3269e-01 & 1.0331 & 9.9066e-02 & 1.0204 &  1.7027e-01 & 1.0167 &  1.9984e-01 & 1.0240 & 6.6419e-02  & 0.9874 \\
0 &  8749 & 2.7038e-01 & 1.0166 & 6.2010e-02 & 1.0129 &  1.0674e-01 & 1.0096 &  1.2490e-01 & 1.0162 & 4.1913e-02  & 0.9954 \\
  & 19605 & 1.8014e-01 & 1.0066 & 4.1346e-02 & 1.0047 &  7.1146e-02 & 1.0056 &  8.3291e-02 & 1.0043 & 2.8040e-02  & 0.9963 \\
  & 43805 & 1.2010e-01 & 1.0086 & 2.7577e-02 & 1.0075 &  4.7506e-02 & 1.0047 &  5.5492e-02 & 1.0103 & 1.8769e-02  & 0.9987 \\ \hline
  &  1766 & 5.8895e-02 &  $--$  & 2.6740e-02 &  $--$  &  2.4564e-02 &  $--$  &  1.7131e-02 &  $--$  & 1.0741e-02  &  $--$  \\
  & 10406 & 9.5028e-03 & 2.0569 & 4.3181e-03 & 2.0560 &  4.0707e-03 & 2.0268 &  2.3959e-03 & 2.2182 & 1.7408e-03  & 2.0519 \\
1 & 26246 & 3.7187e-03 & 2.0283 & 1.6914e-03 & 2.0262 &  1.6028e-03 & 2.0150 &  8.9830e-04 & 2.1208 & 6.8185e-04  & 2.0263 \\
  & 58814 & 1.6535e-03 & 2.0089 & 7.5239e-04 & 2.0079 &  7.1485e-04 & 2.0014 &  3.8855e-04 & 2.0774 & 3.0292e-04  & 2.0111 \\
  & 131414& 7.3548e-04 & 2.0154 & 3.3486e-04 & 2.0139 &  3.1868e-04 & 2.0097 &  1.6997e-04 & 2.0569 & 1.3484e-04  & 2.0135 \\ \hline
  &  3384 & 2.1867e-03 &  $--$  & 1.5283e-03 &  $--$  &  7.9424e-04 &  $--$  & 2.3731e-03 &  $--$  & 4.2607e-04  &  $--$  \\
  & 19944 & 1.4347e-04 & 3.0713 & 9.9799e-05 & 3.0766 &  5.3189e-05 & 3.0482 & 1.5168e-04 & 3.1008 & 2.7317e-05  & 3.0973 \\ 
2 & 50304 & 3.5232e-05 & 3.0356 & 2.4480e-05 & 3.0379 &  1.3119e-05 & 3.0261 & 3.7011e-05 & 3.0494 & 6.6740e-06  & 3.0466 \\
  & 112726& 1.0462e-05 & 3.0096 & 7.2624e-06 & 3.0120 &  3.9066e-06 & 3.0027 & 1.0949e-05 & 3.0189 & 1.9758e-06  & 3.0172 \\
  & 251876& 3.1074e-06 & 3.0200 & 2.1565e-06 & 3.0206 &  1.1618e-06 & 3.0168 & 3.2445e-06 & 3.0257 & 5.8587e-07  & 3.0240 \\ \hline
\end{tabular} }\vspace{0.5cm}   
\caption{ Test 1. Convergence history of the terms composing the estimator using hexagonal meshes. }  
\label{tab2} 
\end{table}

\subsection{Test 2. Solution with a sharp layer: uniform vs adaptive refinement}

\begin{figure}
\begin{center}
\scalebox{0.3}{ \includegraphics{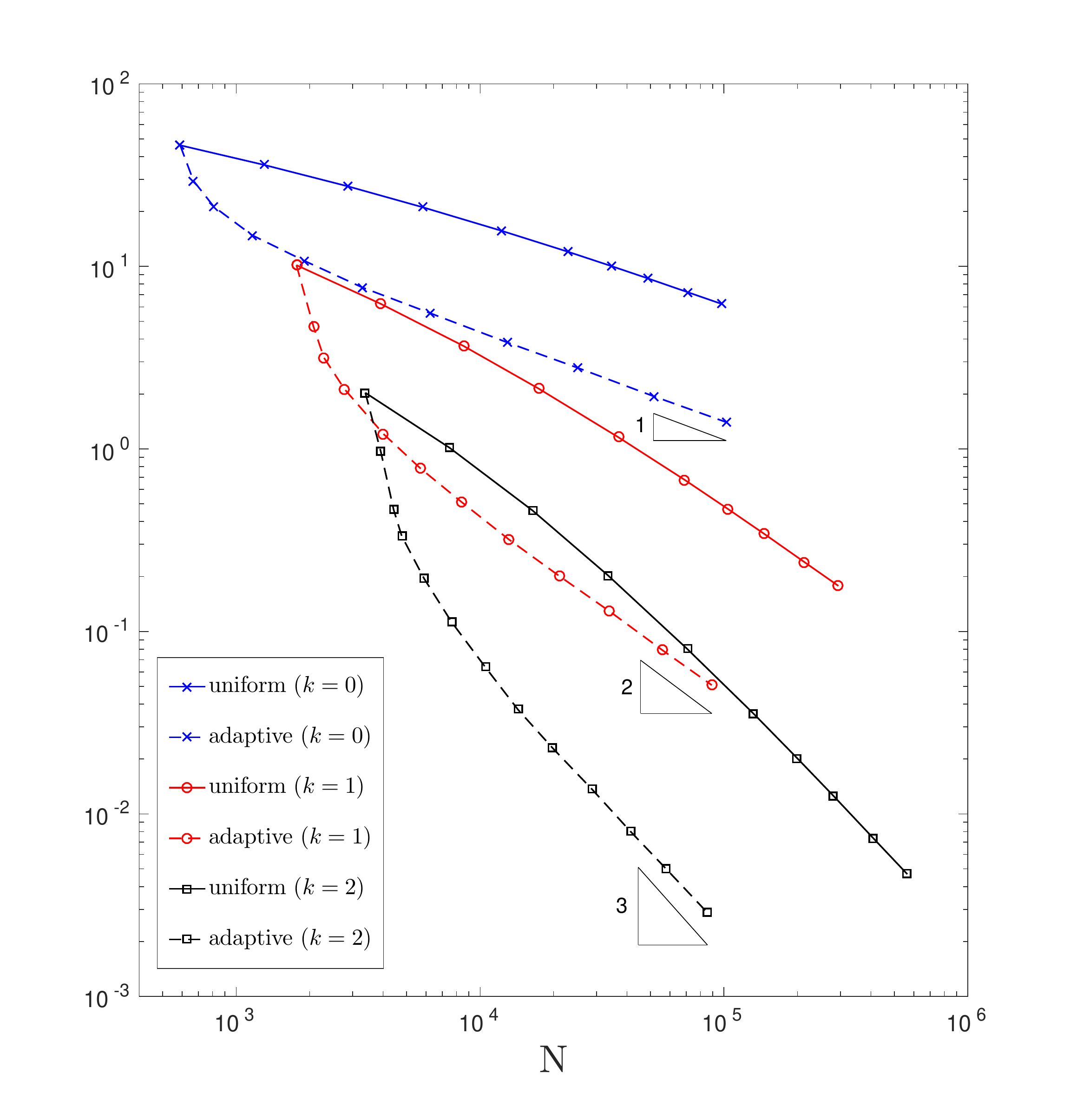}\hspace{-1cm}
                 \includegraphics{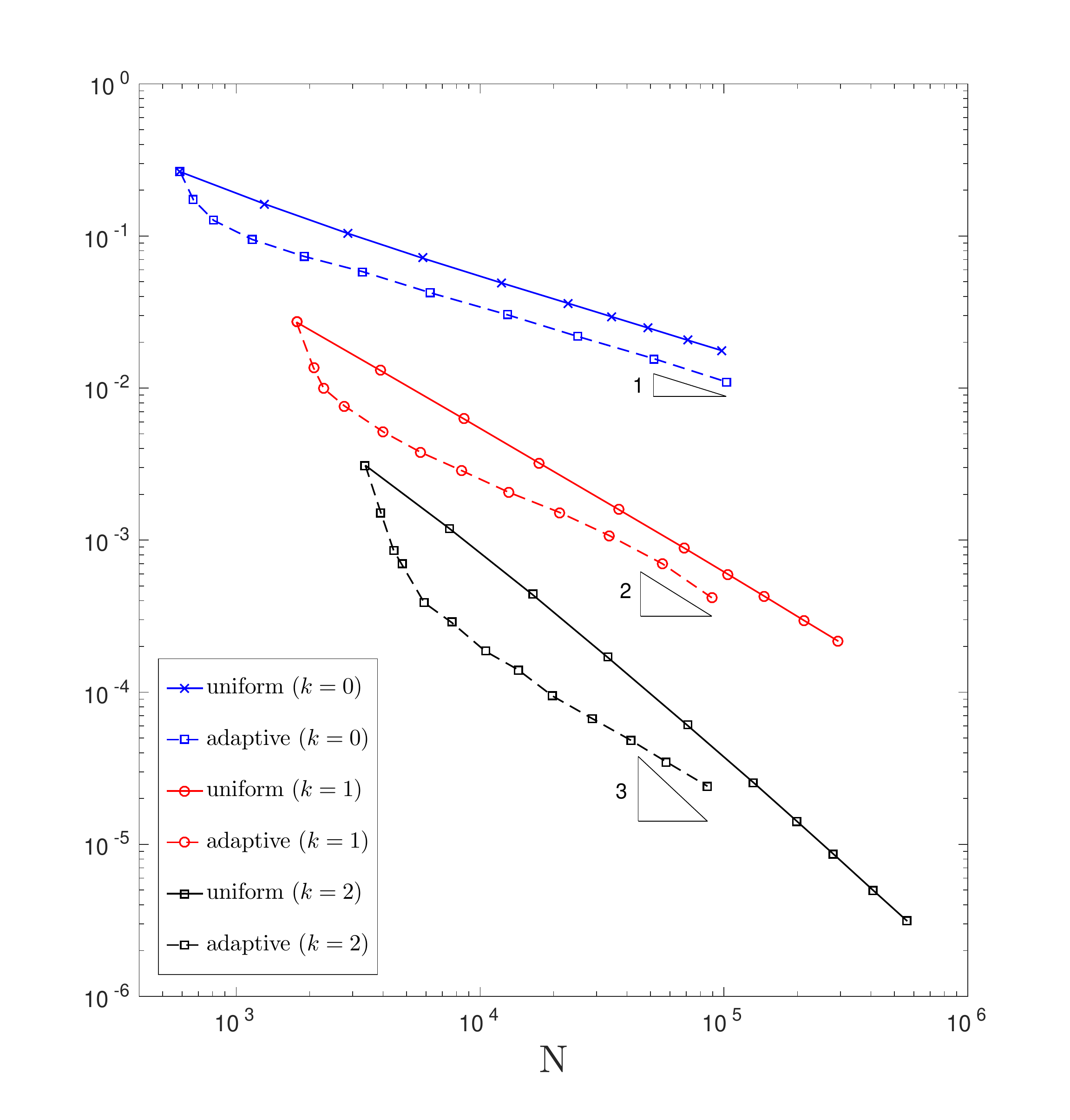}}\vspace{-0.5cm}
\end{center}
\caption{Test 2. Convergence history under uniform  and the adaptive refinement of hexagonal meshes (cf. Figure~\ref{pict2.5}). The error ${\tt e}(\bsi)$ (left) and  ${\tt e}(u)$ (right).}\label{pict2.1}  
\end{figure}

\begin{figure}
\begin{center}
\scalebox{0.25}{
\hspace{-2cm}
\includegraphics{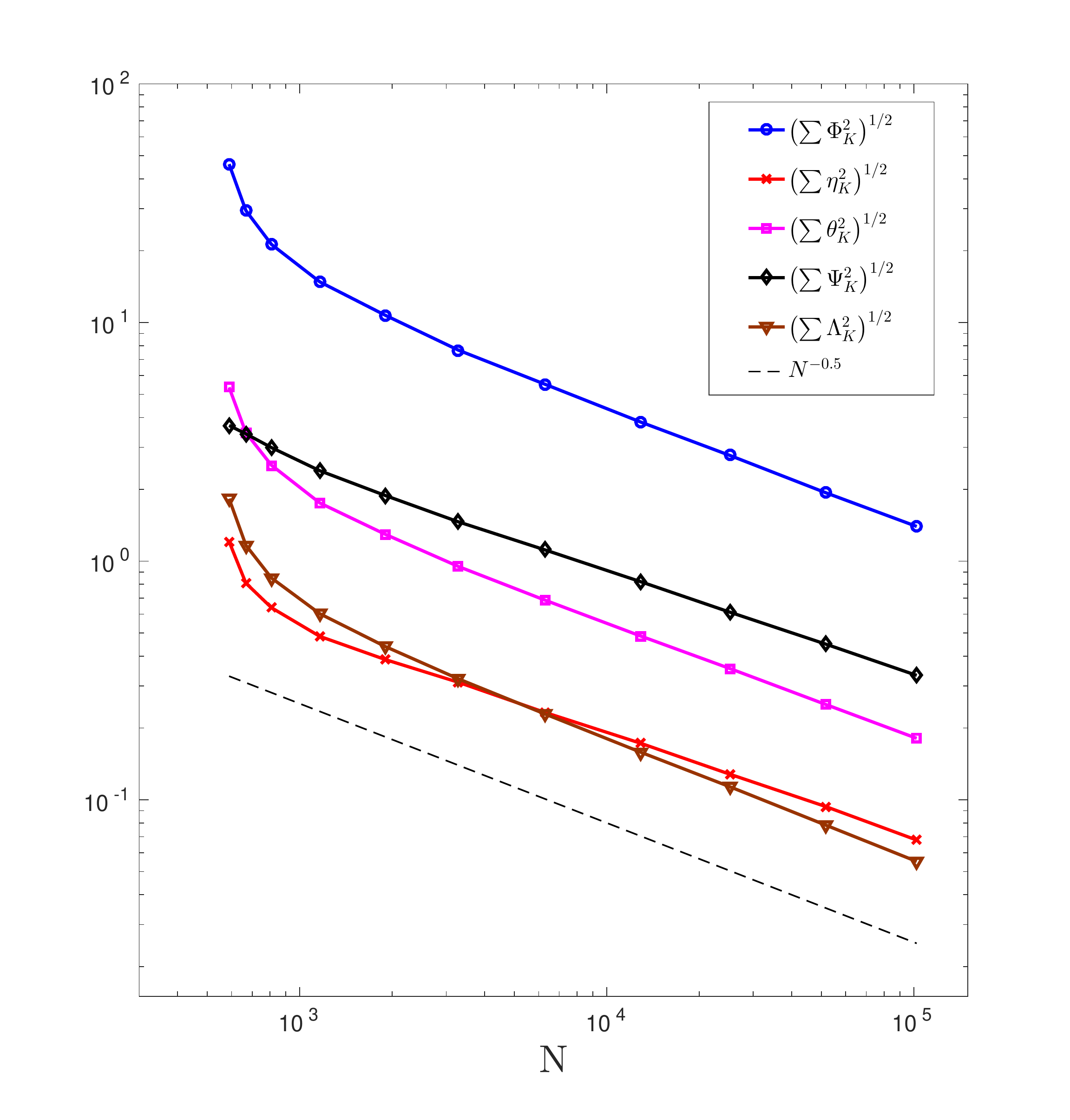}\hspace{-1.5cm}
                 \includegraphics{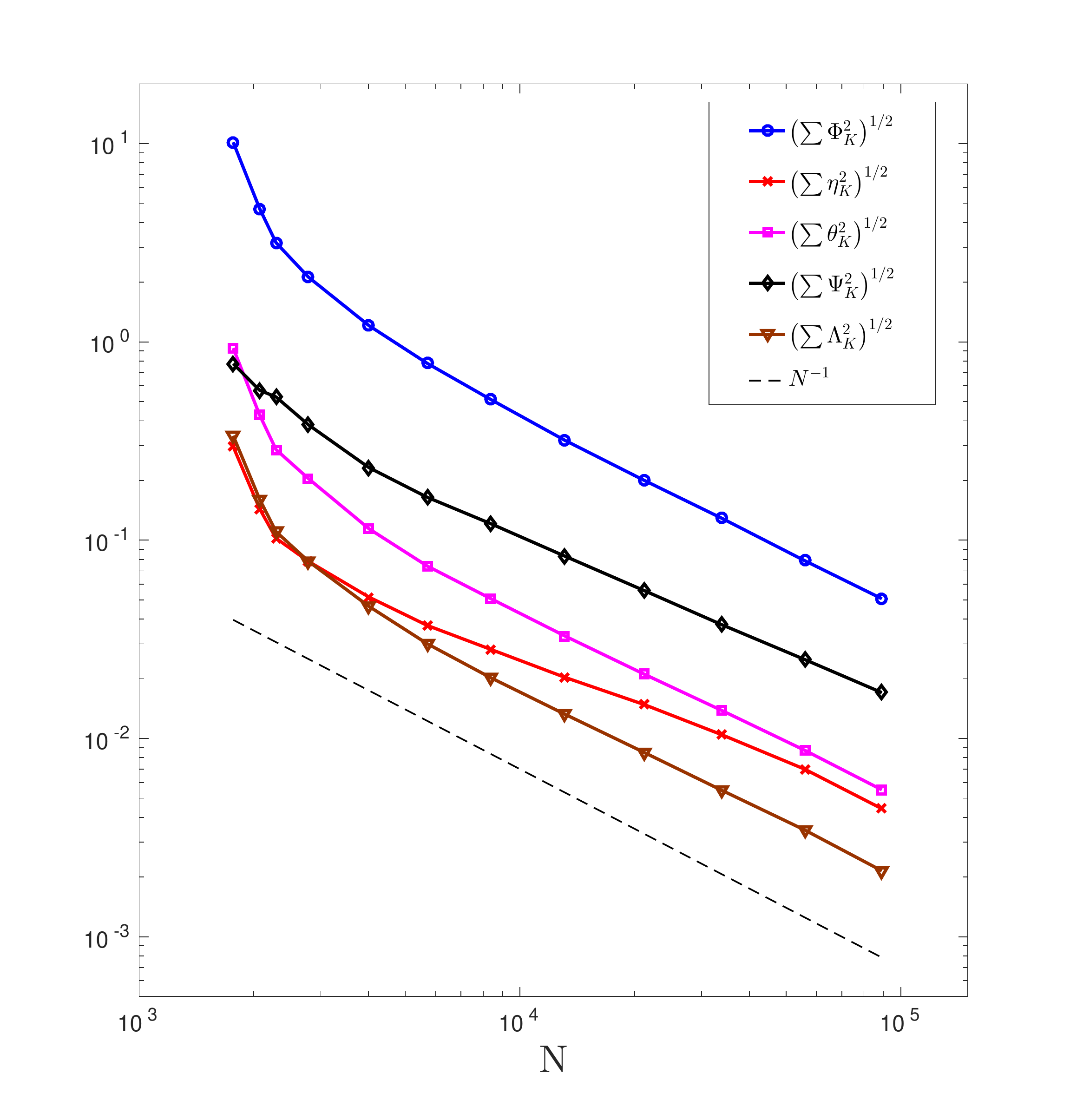}\hspace{-1.5cm}
                 \includegraphics{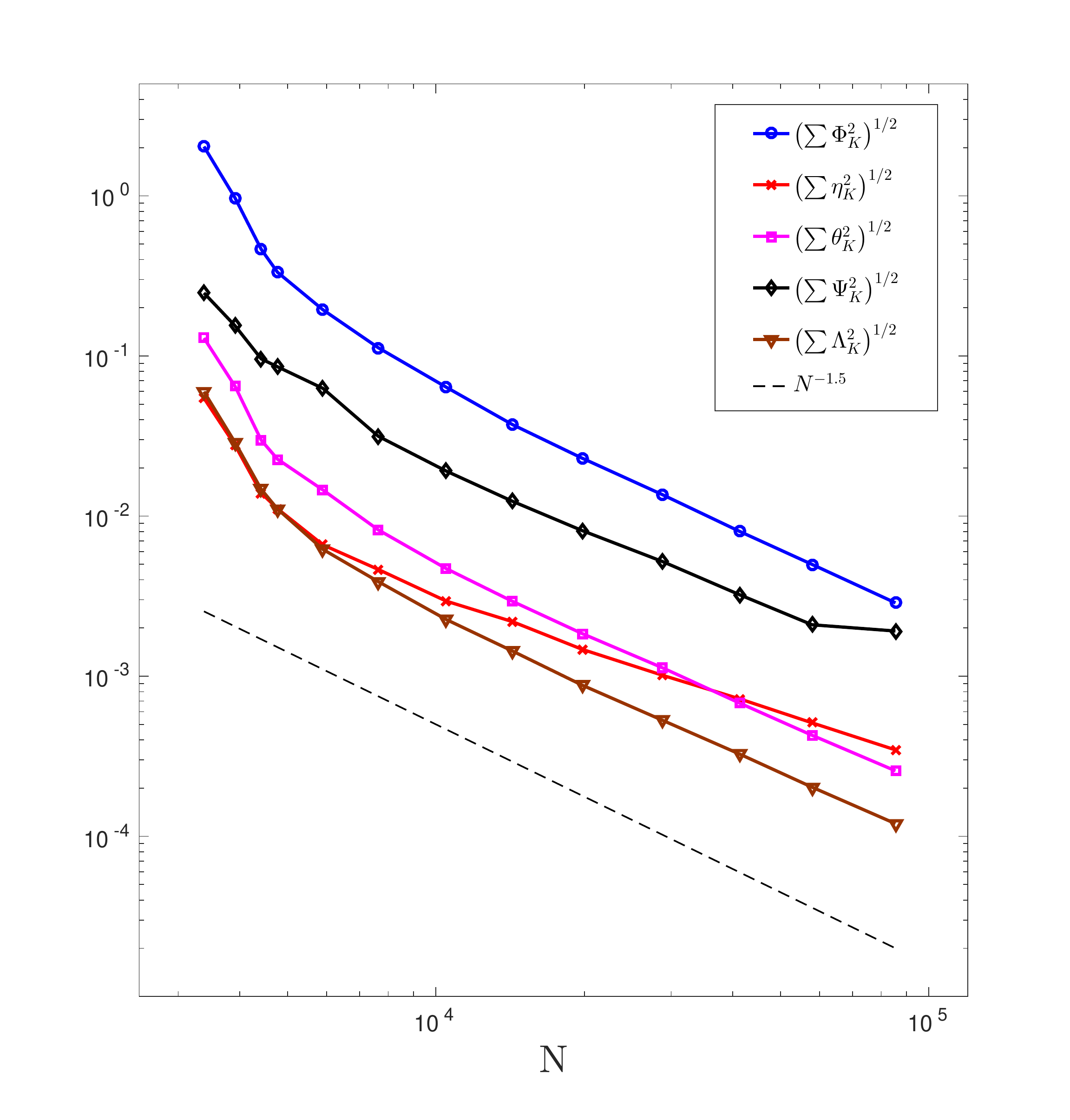} }\vspace{-0.5cm}           
\end{center}
\caption{Test 2. Convergence history of the components of the
estimator under adaptive refinement of hexagonal meshes (cf. Figure~\ref{pict2.5} below).  For $k=0$ (left), $k=1$ (centre), and $k=2$ (right).}\label{pict2.2}  
\end{figure}

We consider $\Omega:=(0,1)^2$ with $\G_D:=\Big\{(w,0),(0,w)\in\Omega:\quad 0\leq w\leq 1 \Big\}$ and $\G_N:=\G\setminus\overline{\G}_D$ , and choose $f$ and $g$ such that the exact solution is given by
\begin{equation*}
u(x,y) =(x-1)^2(y-1)^2\left(\frac{1}{x+0.1} + \frac{1}{1+y}\right)\qin\Omega\,.
\end{equation*}
Note that $u$ and $\nabla u$ are singular along the lines 
$x =\,-\,0.1$ and $y =\,-\,1$. Both such lines are outside $\Omega$, but we expect regions of high gradients in the vicinity of the left boundary.
From Figure \ref{pict2.1} we observe, as expected, that  the adaptive methods outperforms  uniform refinement. Indeed, initially the adaptive method superconverges until, ones the steep layer is resolved, both methods converge at the theoretical rate, namely $k+1$.
This is clearly shown in Table \ref{tab3}, where the rates of convergence of the global error and the estimator at each step of the adaptive algorithm are reported together with the effectivity index.
 As shown in Figure \ref{pict2.2}, all terms in the error estimator follow precisely the same behaviour. 

\begin{table}
\centering
\scalebox{0.93}{\begin{tabular}{|c|r|c c|c c|c|}\hline
$k$& $N$ & ${\tt e}(\bsi,u)$ & ${\tt r}(\bsi,u)$ & $\Theta$ & ${\tt r}(\Theta)$ & ${\tt eff}(\Theta)$  \\ \hline
  &    589 & 4.6114e+01 &  $--$  & 4.6592e+01 &  $--$  & 0.9896  \\
  &    668 & 2.9471e+01 & 7.1146 & 2.9879e+01 & 7.0595 & 0.9863  \\
  &    809 & 2.1279e+01 & 3.4013 & 2.1642e+01 & 3.3681 & 0.9832  \\
0 &   1163 & 1.4820e+01 & 1.9932 & 1.5118e+01 & 1.9769 & 0.9803  \\
  &   1902 & 1.0735e+01 & 1.3112 & 1.0978e+01 & 1.3011 & 0.9779  \\
  &   3290 & 7.6690e+00 & 1.2275 & 7.8661e+00 & 1.2165 & 0.9749 \\
  &   6272 & 5.5309e+00 & 1.0131 & 5.6850e+00 & 1.0066 & 0.9729  \\
  &  12928 & 3.8304e+00 & 1.0158 & 3.9474e+00 & 1.0086 & 0.9704  \\
\hline  
  &   1766 & 1.0162e+01 &  $--$  & 1.0241e+01 &  $--$  & 0.9923  \\
  &   2072 & 4.7026e+00 & 9.6441 & 4.7595e+00 & 9.5903 & 0.9880  \\
  &   2288 & 3.1654e+00 & 7.9834 & 3.2240e+00 & 7.8561 & 0.9818  \\
  &   2782 & 2.1237e+00 & 4.0834 & 2.1691e+00 & 4.0545 & 0.9791  \\
1 &   4014 & 1.2081e+00 & 3.0775 & 1.2367e+00 & 3.0652 & 0.9769  \\
  &   5706 & 7.8010e-01 & 2.4868 & 8.0163e-01 & 2.4652 & 0.9732  \\
  &   8368 & 5.1334e-01 & 2.1859 & 5.3063e-01 & 2.1550 & 0.9674  \\
  &  13090 & 3.1982e-01 & 2.1151 & 3.3270e-01 & 2.0867 & 0.9613  \\
  &  21158 & 2.0077e-01 & 1.9394 & 2.0991e-01 & 1.9183 & 0.9564  \\
\hline
  &   3384 & 2.0312e+00 &  $--$  & 2.0513e+00 &  $--$  & 0.9902  \\
  &   3913 & 9.7026e-01 &10.1734 & 9.8510e-01 &10.0996 & 0.9849  \\
  &   4422 & 4.6358e-01 &12.0796 & 4.7458e-01 &11.9440 & 0.9768  \\
  &   4771 & 3.3337e-01 & 8.6810 & 3.4515e-01 & 8.3847 & 0.9659  \\
2 &   5899 & 1.9436e-01 & 5.0845 & 2.0487e-01 & 4.9154 & 0.9487  \\
  &   7640 & 1.1258e-01 & 4.2231 & 1.1729e-01 & 4.3134 & 0.9598  \\
  &  10494 & 6.3836e-02 & 3.5747 & 6.6868e-02 & 3.5407 & 0.9547  \\
  &  14293 & 3.7310e-02 & 3.4766 & 3.9482e-02 & 3.4105 & 0.9450  \\
  &  19800 & 2.3005e-02 & 2.9673 & 2.4509e-02 & 2.9259 & 0.9386  \\
\hline
\end{tabular} } \vspace{0.5cm}  
\caption{Test 2. The behaviour of the global error and the estimator under adaptive refinement of hexagonal meshes (cf. Figure~\ref{pict2.5}). The effectivity of the estimator is reported in the right-most column.}  
\label{tab3}  
\end{table}

\begin{figure}
\begin{center}
\scalebox{0.25}{ \includegraphics{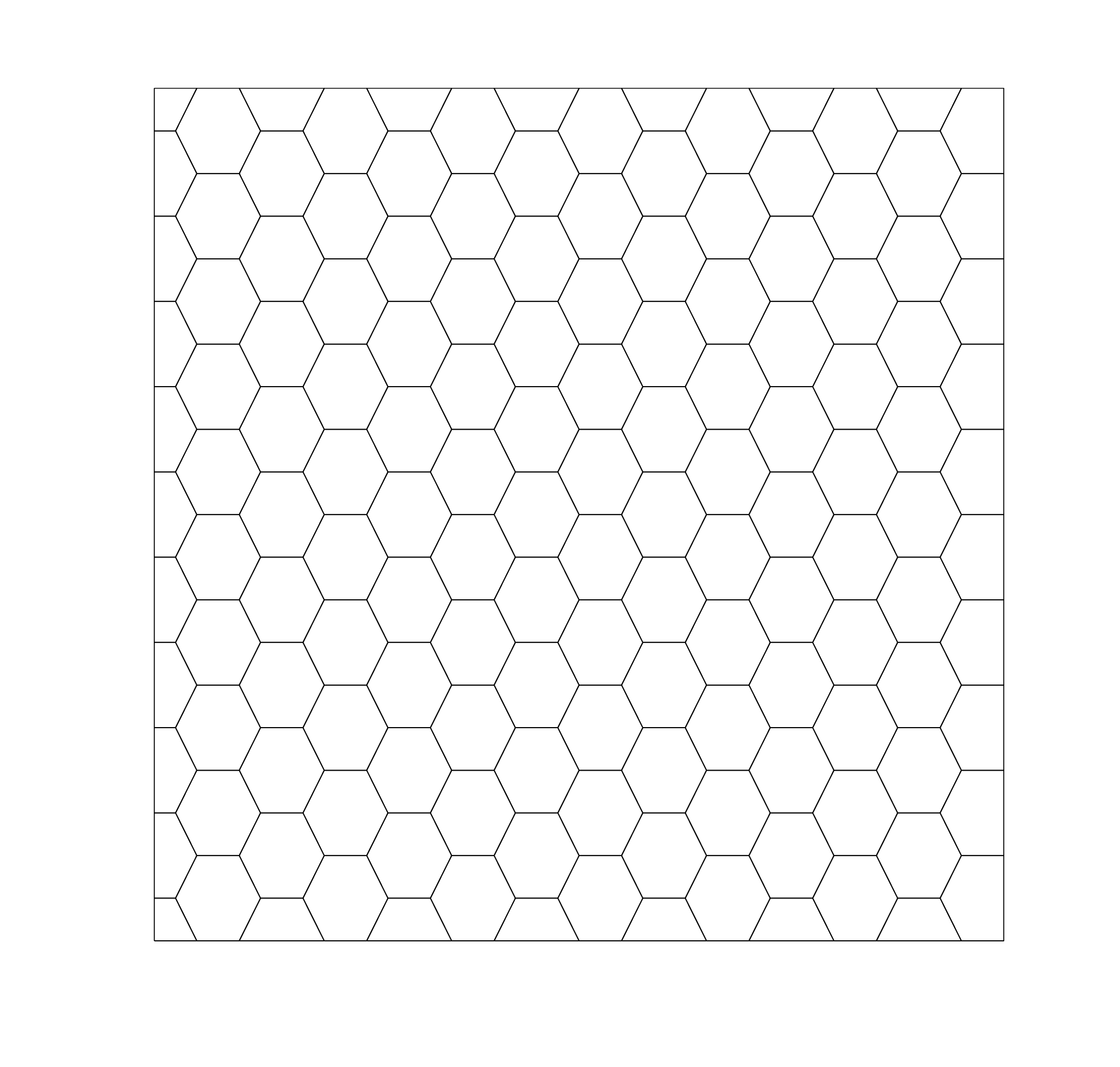}\hspace{-2cm}
                \includegraphics{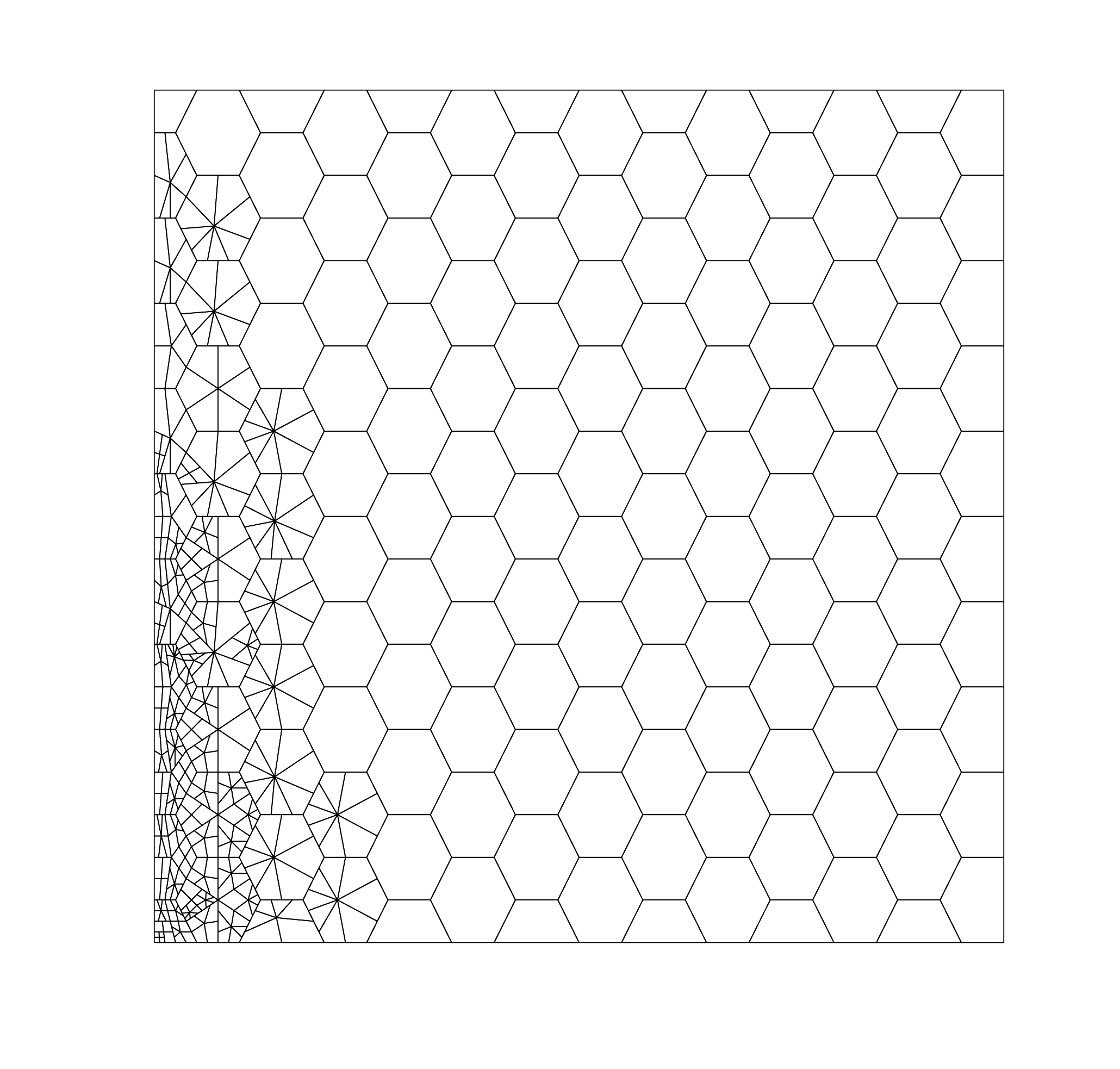}\hspace{-2cm}
                \includegraphics{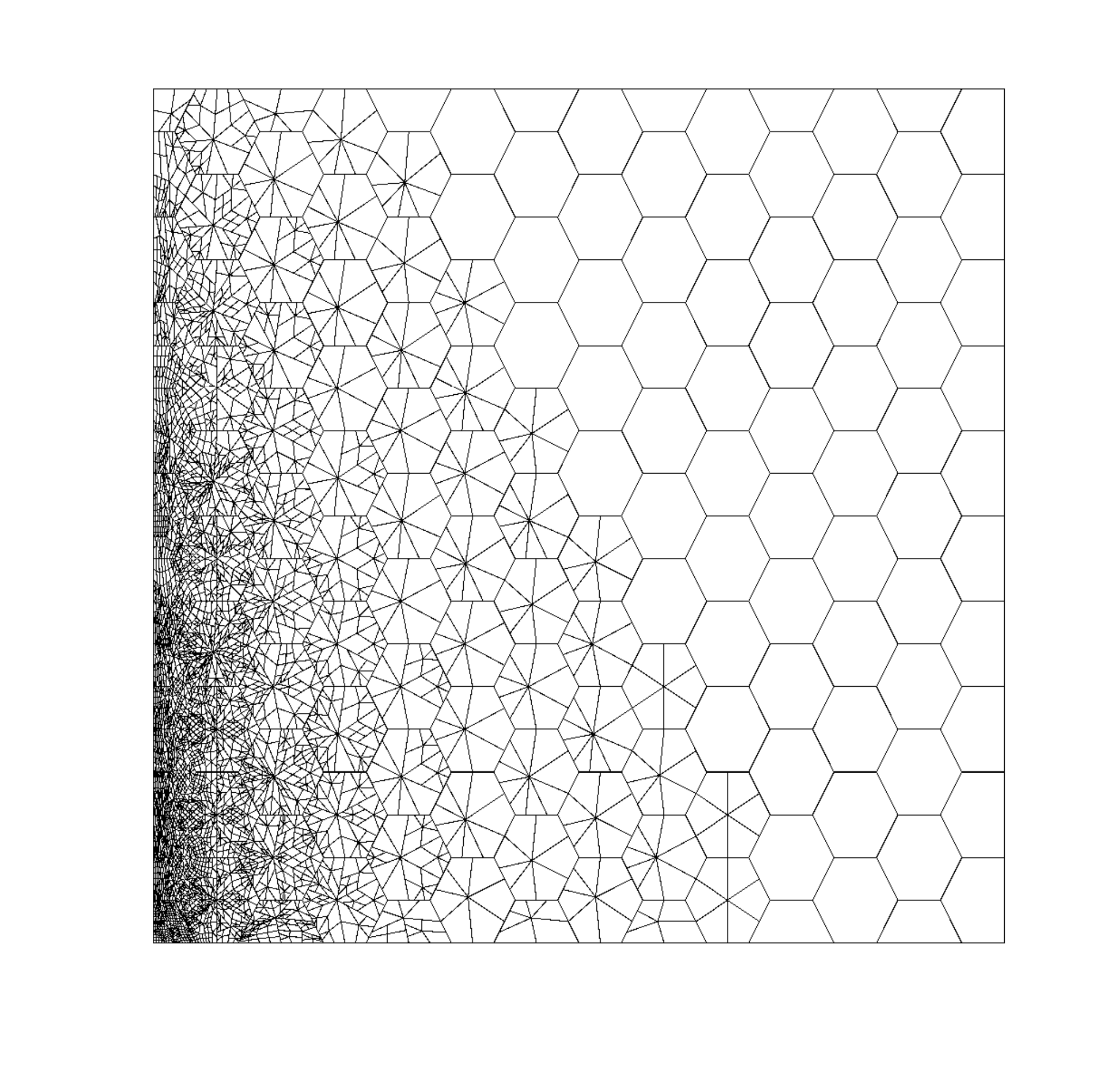}}\vspace{-0.5cm}  
\end{center}
\caption{Test 2. Some meshes from the adaptive refinement sequence obtained with $k=1$: initial (left), after $5$ refinement steps (centre), and after $10$ refinement steps (right). }\label{pict2.5}  
\end{figure}

Some intermediate meshes
obtained with adaptive strategy are displayed in Figure~\ref{pict2.5}. Notice here that the adapted meshes concentrate the refinements in the proximity of  the line $x= 0$, confirming that the adaptive algorithm is able to target the regions with high gradients of the solution.

\subsection{Test 3. L-shaped domain solution: adaptive refinement}

\begin{figure}[h]
\begin{center}
\scalebox{0.35}{
\includegraphics{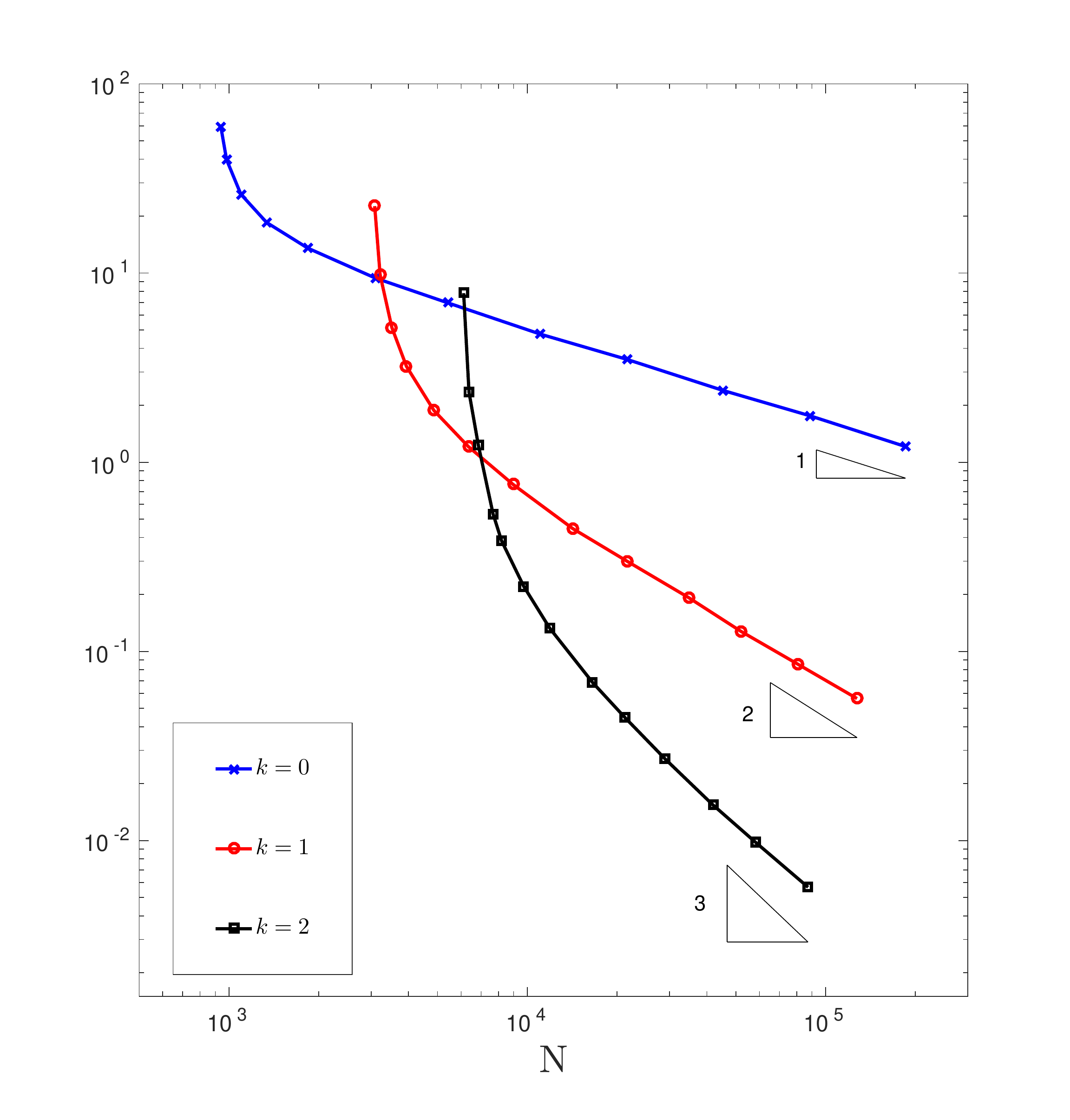}
\includegraphics{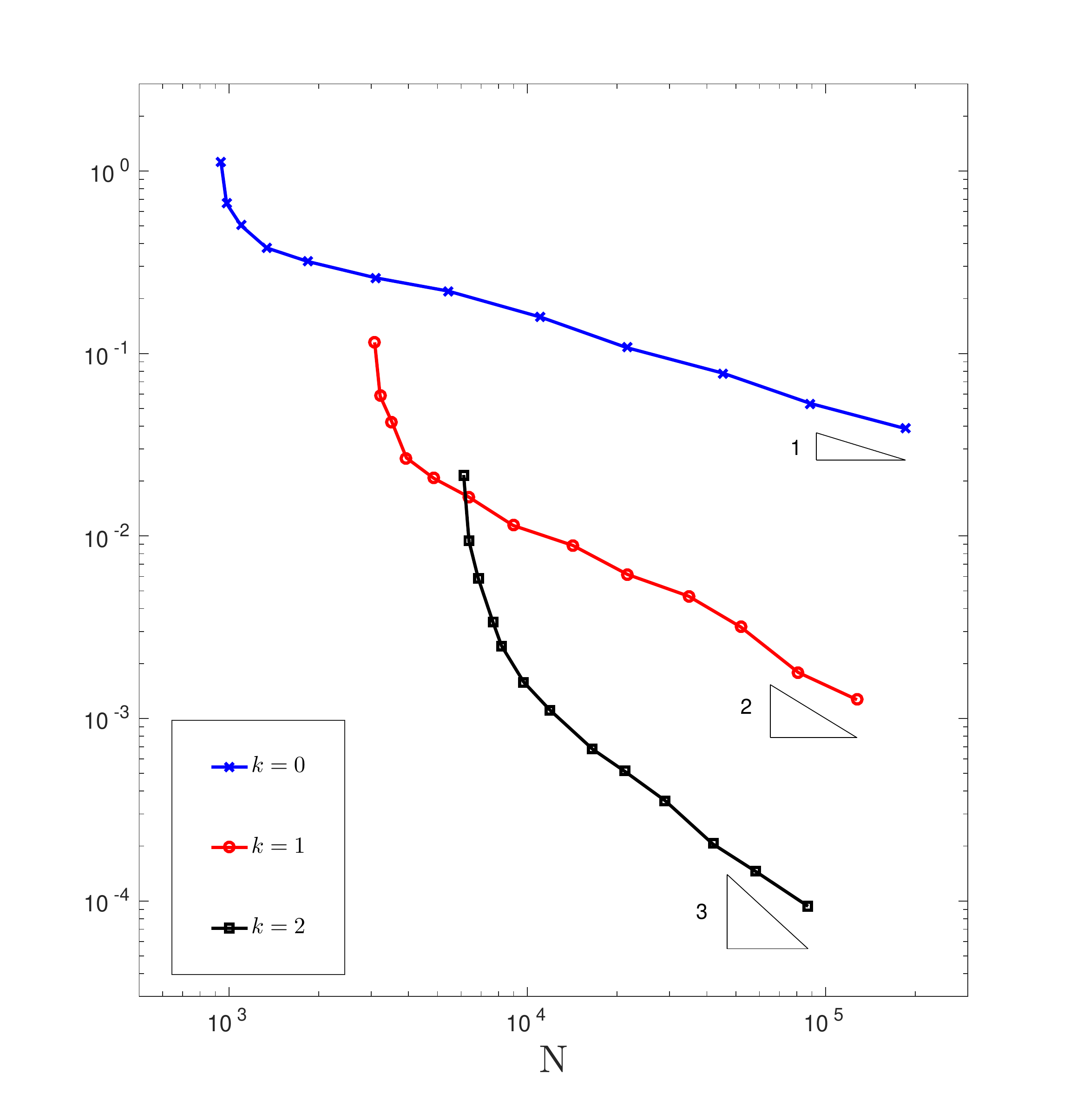}}\vspace{-0.5cm}
\end{center}
\caption{Test 3. Errors curves for the adaptive strategy using distorted quadrilateral meshes, (cf. Figure~\ref{pict3.3} below). The error ${\tt e}(\bsi)$ (left) and the error ${\tt e}(u)$ (right). }\label{pict3.1}  
\end{figure}

\begin{figure}
\begin{center}
\includegraphics[width=.31\textwidth]{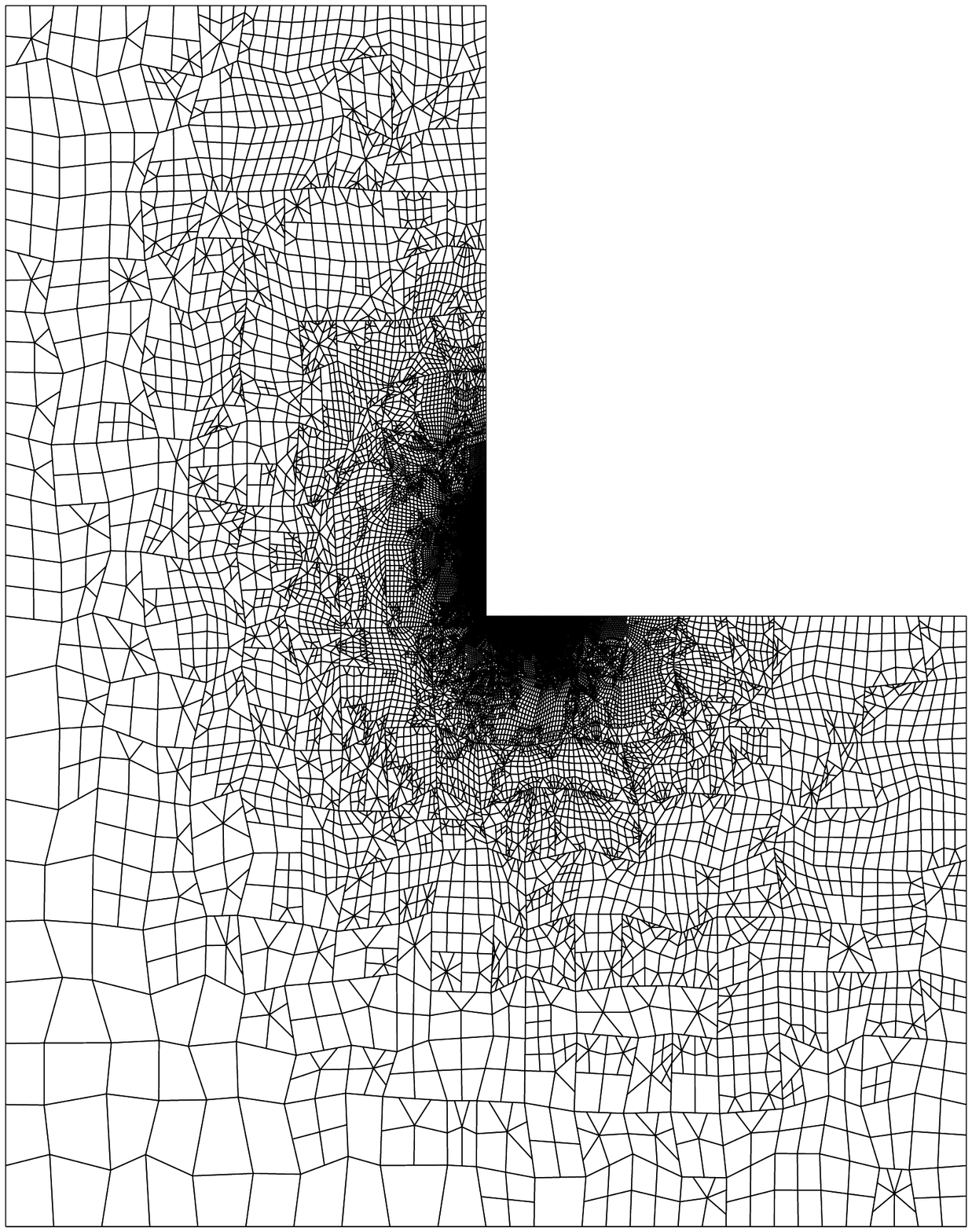}\hspace{-.4cm}
                \includegraphics[width=.31\textwidth]{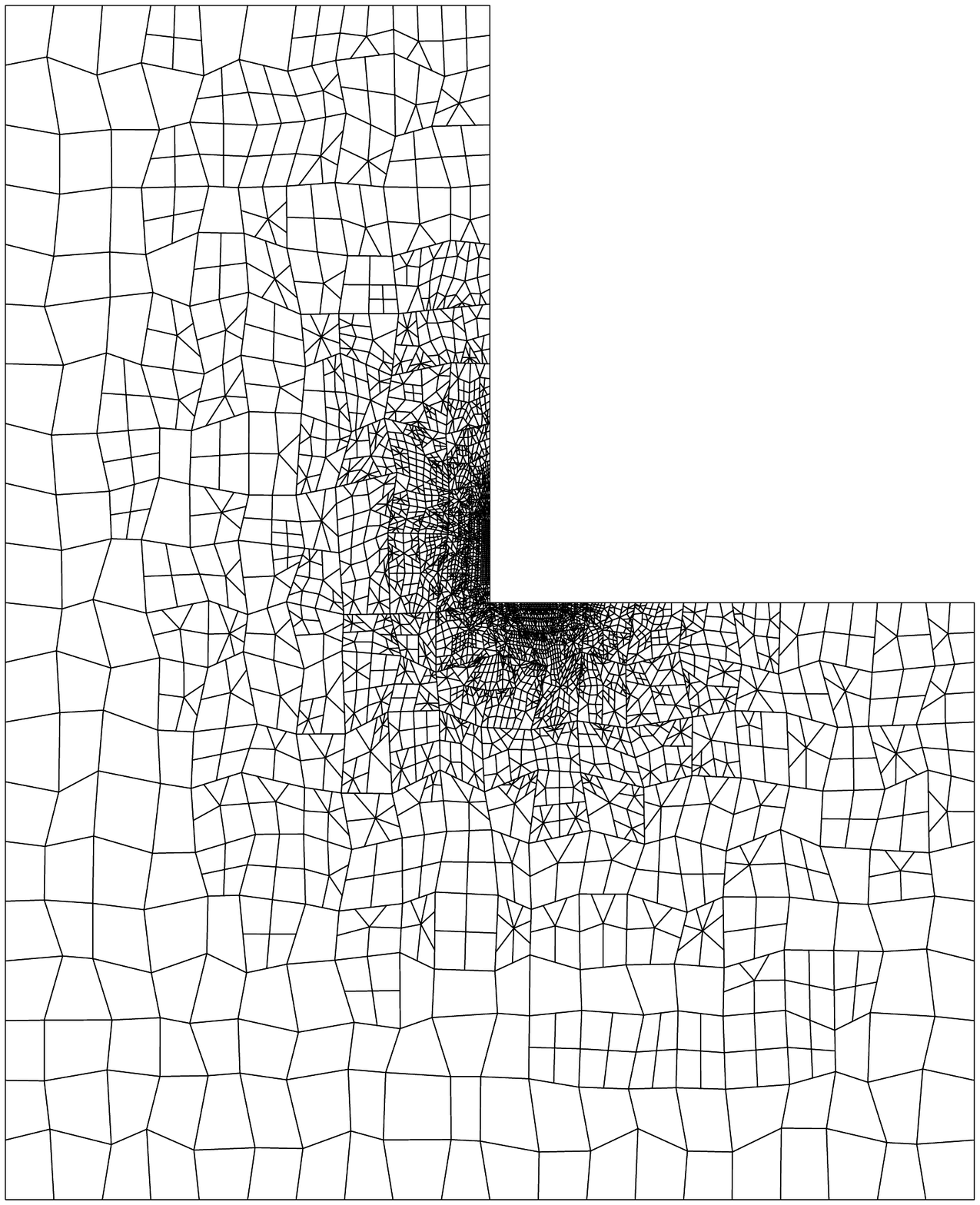}\hspace{-.4cm}
                \includegraphics[width=.31\textwidth]{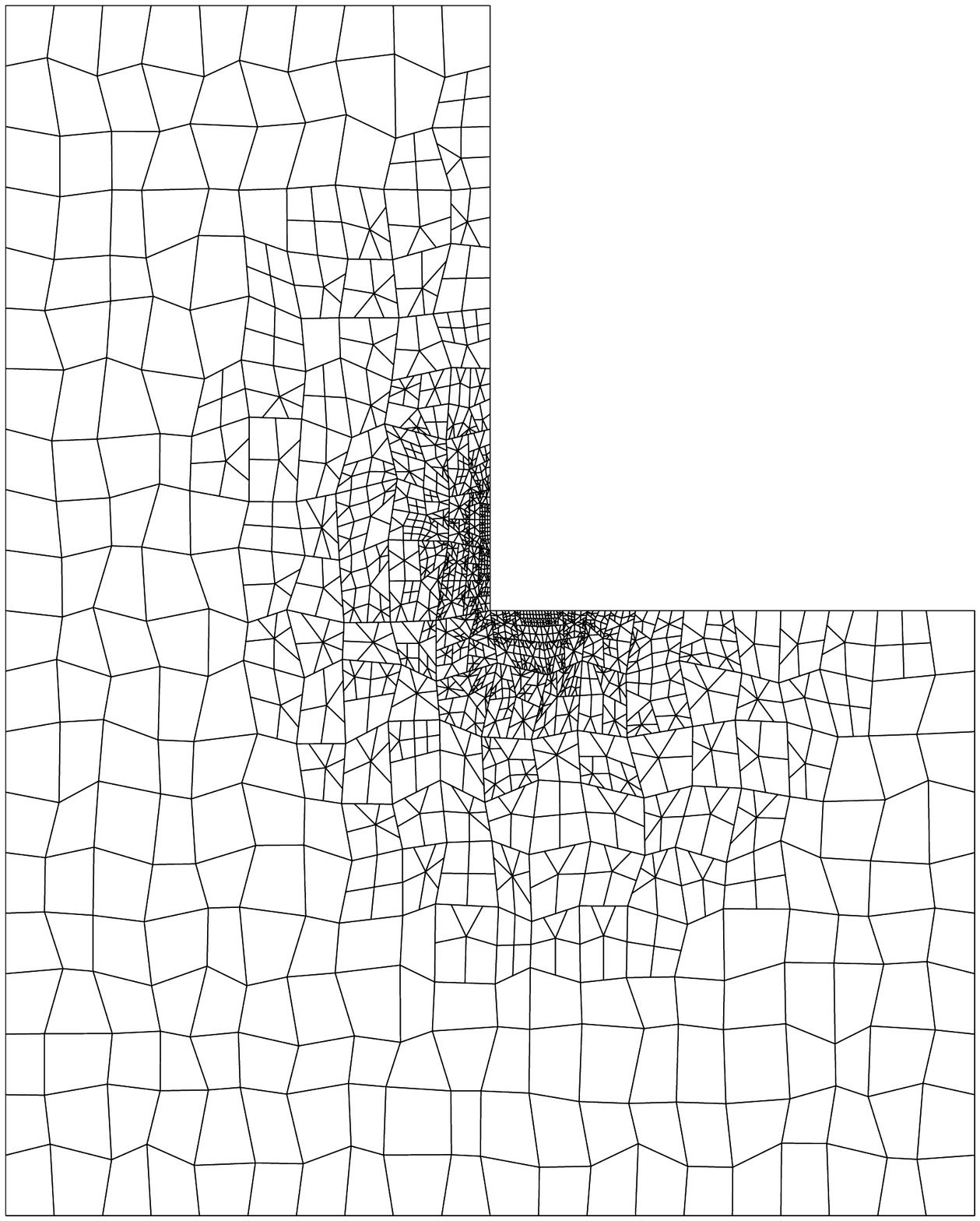}\vspace{-1cm}
\includegraphics[width=.31\textwidth]{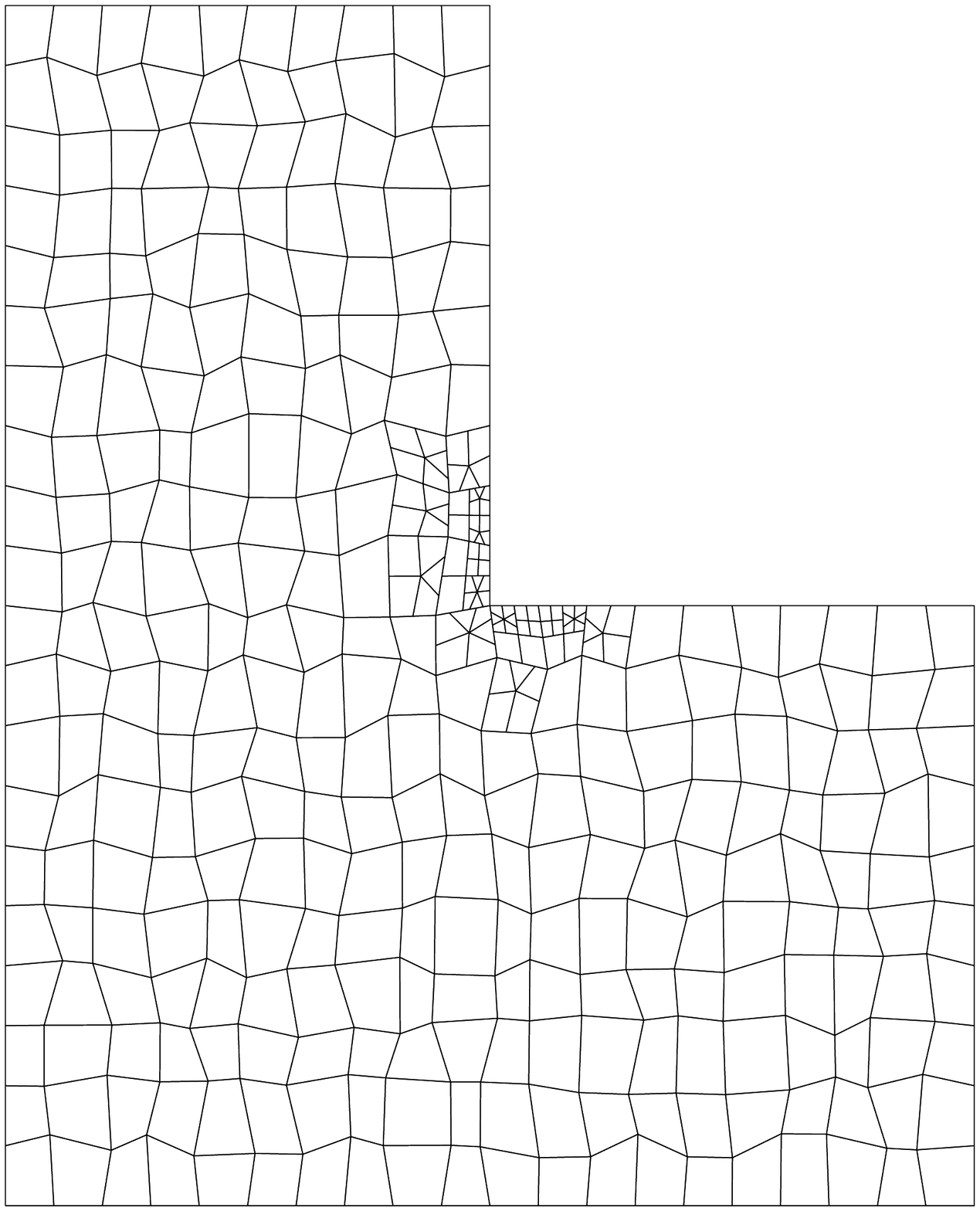}\hspace{-.5cm}
                \includegraphics[width=.31\textwidth]{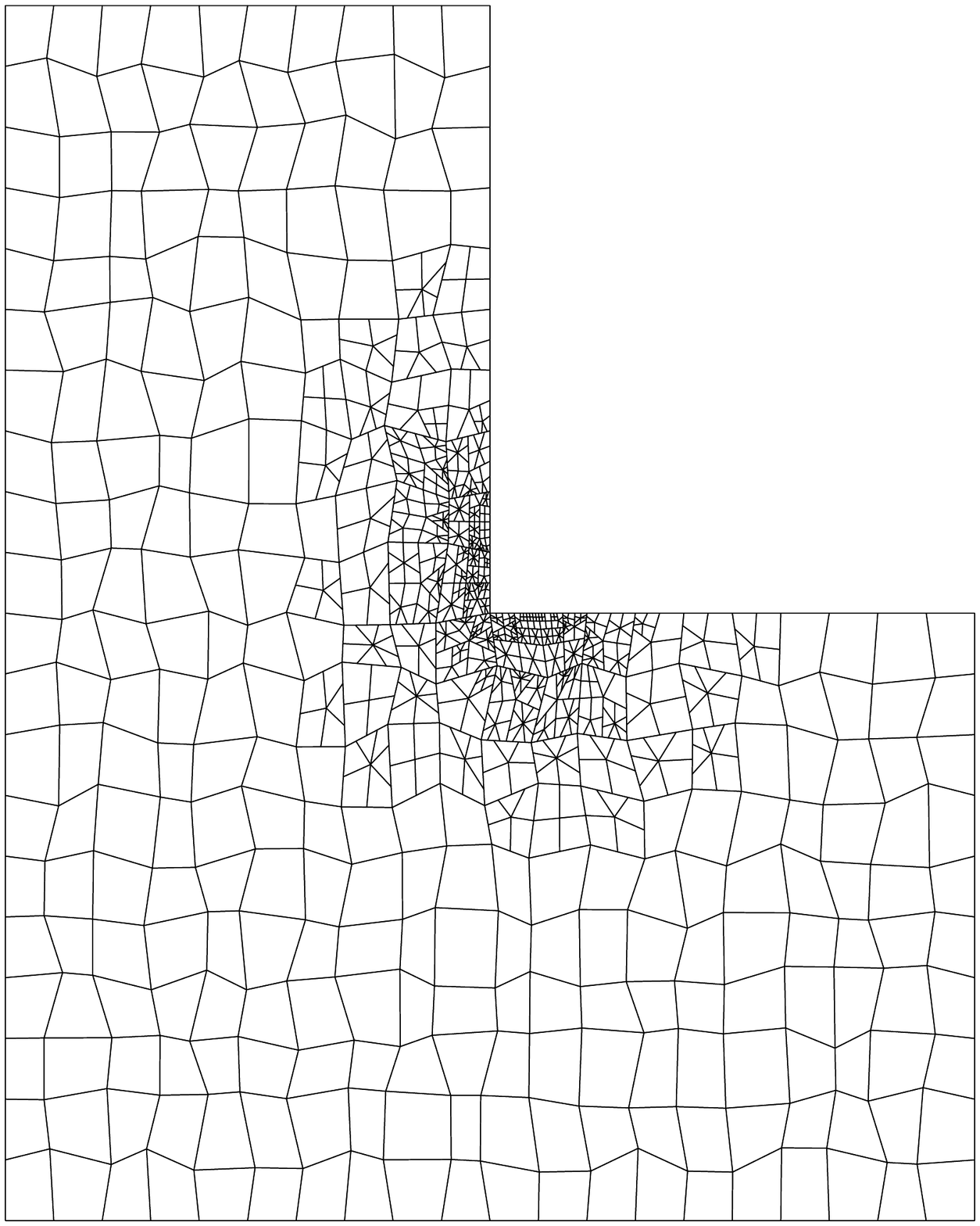}\hspace{-.5cm}
                \includegraphics[width=.31\textwidth]{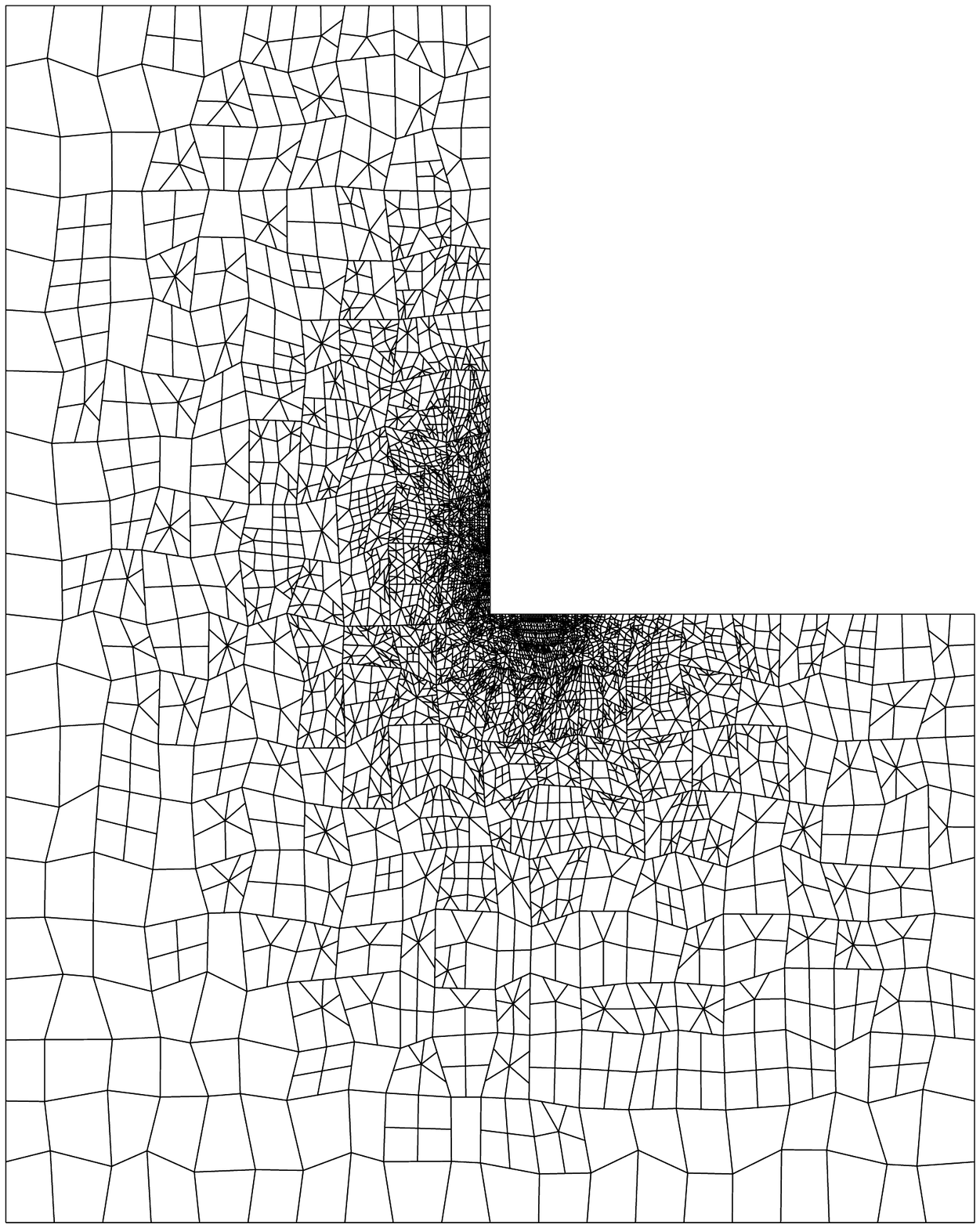}
\vspace{-1cm}         
\end{center}
\caption{Test 3. (Top)  The mesh after ten adaptive refinements with $k=0$ (left), $k=1$ (centre) and $k=2$ (right). (Below) Some meshes from the adaptive refinement sequence for $k=2$: after $3$ refinement steps (left), $8$ refinement steps (centre), and $15$ refinement steps (right).}\label{pict3.3}  
\end{figure}

We consider $\Omega:=(-1,-1)^2\setminus (0,1)^2$ with $\G_N:=\Big\{(-1,w),(w,-1)\in\Omega:\quad -1\leq w\leq 1 \Big\}$ and $\G_D:=\G\setminus\overline{\G}_N$,  and choose $f$ and $g$ such that the exact solution is given by
\begin{equation*}
u(x,y) =\dfrac{(x+1)^2(y+1)^2}{\sqrt{(x-0.1)^2 + (y-0.1)^2}}\qin\Om\,.
\end{equation*}
Note that $\Omega$ is an L-shaped domain and that $u$ and $\nabla u$  are singular at the point $(0.1, 0.1)$, which is just outside of $\Omega$. Hence, we should expect regions of high gradients around the origin, which
is the middle corner of the L-shaped domain.  
In Figure \ref{pict3.1} and Table \ref{tab4} we display the convergence history of the
adaptive method. Finally, Figure~\ref{pict3.3} shows how the adaptive strategy correctly refines in a neighbourhood of the origin. We also notice that increasing the order of the method  allows for a less aggressive refinement.

\begin{table}
\centering
\scalebox{0.93}{\begin{tabular}{|c|r|c c|c c|c|}\hline
$k$& $N$ & ${\tt e}(\bsi,u)$ & ${\tt r}(\bsi,u)$ & $\Theta$ & ${\tt r}(\Theta)$ & ${\tt eff}(\Theta)$  \\ \hline
  &    940 & 5.9057e+01 &  $--$  & 6.1812e+01 &  $--$  & 0.9554  \\
  &    982 & 3.9730e+01 &18.1365 & 4.1712e+01 &17.9953 & 0.9525  \\
  &   1096 & 2.6136e+01 & 7.6263 & 2.7724e+01 & 7.4387 & 0.9427  \\
0 &   1337 & 1.8521e+01 & 3.4653 & 2.0096e+01 & 3.2376 & 0.9216  \\
  &   1838 & 1.3548e+01 & 1.9650 & 1.4905e+01 & 1.8783 & 0.9090  \\
  &   3098 & 9.4108e+00 & 1.3959 & 1.0536e+01 & 1.3286 & 0.8932 \\
  &   5420 & 6.9603e+00 & 1.0786 & 7.8716e+00 & 1.0426 & 0.8842  \\
  &  11100 & 4.7501e+00 & 1.0659 & 5.4231e+00 & 1.0395 & 0.8759  \\
\hline  
  &   3080 & 2.2588e+01 &  $--$  & 2.4078e+01 & $--$   & 0.9381 \\
  &   3212 & 9.7746e+00 &39.9209 & 1.0379e+01 &40.1065 & 0.9418 \\
  &   3516 & 5.1580e+00 &14.1379 & 5.6331e+00 &13.5160 & 0.9157 \\
  &   3936 & 3.2117e+00 & 8.3964 & 3.5285e+00 & 8.2911 & 0.9102 \\
1 &   4866 & 1.8770e+00 & 5.0650 & 2.1510e+00 & 4.6667 & 0.8726 \\
  &   6342 & 1.2110e+00 & 3.3083 & 1.4262e+00 & 3.1021 & 0.8491 \\
  &   8966 & 7.6293e-01 & 2.6687 & 9.0810e-01 & 2.6076 & 0.8401 \\
  &  14202 & 4.4611e-01 & 2.3334 & 5.5115e-01 & 2.1714 & 0.8094 \\
  &  21684 & 2.9875e-01 & 1.8948 & 3.6841e-01 & 1.9037 & 0.8109 \\
\hline 
  &   6120 & 7.8275e+00 &  $--$  & 8.6432e+00 & $--$   & 0.9056 \\
  &   6378 & 2.3483e+00 &58.3131 & 2.5850e+00 &58.4627 & 0.9084 \\
  &   6851 & 1.2364e+00 &17.9329 & 1.4322e+00 &16.5094 & 0.8633 \\
  &   7676 & 5.3466e-01 &14.7464 & 6.6401e-01 &13.5204 & 0.8052 \\
2 &   8200 & 3.8435e-01 & 9.9969 & 4.6526e-01 &10.7732 & 0.8261 \\
  &   9738 & 2.2031e-01 & 6.4748 & 2.7509e-01 & 6.1138 & 0.8009 \\
  &  11895 & 1.3259e-01 & 5.0756 & 1.7407e-01 & 4.5744 & 0.7617 \\
  &  16581 & 6.8309e-02 & 3.9936 & 9.0940e-02 & 3.9096 & 0.7512 \\
  &  21183 & 4.4949e-02 & 3.4173 & 6.1021e-02 & 3.2578 & 0.7366 \\
\hline  
\end{tabular} }\vspace{0.5cm}  
\caption{Test 3. The behaviour of the global error and the estimator using the adaptive strategy.  The effectivity of the estimator is reported in the right-most column.}  
\label{tab4}  
\end{table}

\section{Conclusions}\label{sec:conclusions}

We have derived a posteriori error estimates for a mixed-VEM approach for a second order elliptic equation in divergence form with mixed boundary conditions. We have proved upper and lower bounds for the error between the true solution and both the VEM approximation and a computable postprocessing of the VEM approximation. 
In particular, the postprocessing permitted us to obtain optimal error estimates in the broken $\rmH(\div)$-norm, whereas for the directly computable projection of the virtual element approximation, it  is only possible to prove error estimates in the $\rmL^2$-norm. Arguments based in the inf-sup global condition, suitable  Helmholtz  decompositions and a type Cl\'ement-type interpolant were used  to  derive the upper bound.   The lower bound was obtained, in classical fashion,  by using localisation techniques of bubble functions.  We have also proposed an adaptive algorithm based on the fully local and computable error estimator derived from the a posteriori error analysis. Its performance and effectiveness was illustrated through some numerical test. 
The extension of the present analysis to other relevant problems, such as the Stokes system, will be the subject of future works. 

\vspace{0.2cm}
{\it Acknowledgements. } This research was initiated during the visit of AC in Conception, Chile, during a study leave granted by the College of Science and Engineering at the University of Leicester in 2017, and funded by CONICYT-Chile grant No. 1140791, a Santander Travel Grant, and an LMS Caring Suplementary Grant. MM was supported by the Becas-CONICYT Programme for foreign students. All this support is gratefully
acknowledged.

\end{document}